\documentclass[11pt,twoside]{article}
\usepackage{amssymb}
\usepackage{times}
\usepackage{amsmath}
\usepackage[pdftex]{graphicx}
\usepackage{hyperref}

\newtheorem{Thm}{Theorem}[section]
\newtheorem{Thm*}{Theorem}
\newtheorem{Conj}[Thm]{Conjecture}
\newtheorem{Lem}[Thm]{Lemma}

\newtheorem{Def}[Thm]{Definition} 
\newtheorem{Def*}{Definition} 

\newenvironment{Pf}[1]
{\trivlist\item[]{\it #1\@. }}{\hspace*{\fill}$\Box$\endtrivlist}

\renewcommand{\marginpar}[1]{}
\catcode`\@=12

\def\Empty{}
\newcommand\oplabel[1]{
  \def\OpArg{#1} \ifx \OpArg\Empty {} \else
  	\label{#1}
  \fi}
		
\long\def\realfig#1#2#3#4{
\begin{figure}[htbp]
  \caption[#1]{#3}
  \centering
  \includegraphics[width=#4]{#2}
  \label{#1}
\end{figure}}

\newcommand{\comm}[1]{}

\renewcommand{\epsilon}{\varepsilon}

\renewcommand{\rho}{\varrho}

\begin{document}
\title{Monotone homotopies and contracting discs on Riemannian surfaces}
\author{Gregory R. Chambers and Regina Rotman}
\date{October 4, 2016}
\maketitle


\begin{abstract}
	We prove a ``gluing" theorem for monotone homotopies; 
	a monotone homotopy is a homotopy through simple contractible closed curves
	which themselves are pairwise disjoint.
	We show that two monotone homotopies which have appropriate overlap can be replaced
	by a single monotone homotopy.  The ideas used to prove this theorem
	are used in \cite{CL2} to prove an analogous result for cycles,
	which forms a critical step in their proof of the existence of minimal surfaces
	in complete non-compact manifolds of finite volume.

	We also show that, if monotone homotopies exist, then fixed point contractions
	through short curves exist.  In particular, suppose that $\gamma$ is a simple closed
	curve of a Riemannian surface, and that there exists a monotone contraction
	which covers a disc which $\gamma$ bounds consisting of curves of length $\leq L$.
	If $\epsilon > 0$ and $q \in \gamma$, then there exists a homotopy that contracts
	$\gamma$ to $q$ over loops that are based at $q$ and have length bounded by
	$3L + 2d + \epsilon$, where $d$ is the diameter of the surface.
	If the surface is a disc, and if $\gamma$ is the boundary of this disc,
	then this bound can be improved to $L + 2d + \epsilon$.
\end{abstract}

\section*{Introduction.}

The central object of study is a \emph{monotone homotopy} on a Riemannian surface $(M,g)$.
A monotone homotopy is one in which every curve is simple, and in which pairs of curves do not
self-intersect.  For the purposes of this article, we will assume that all curves are
contractible.  If the manifold is a $2$-disc, then this can be rephrased as follows.

\begin{Def} \label{Monotone_disc}
Let $M$ be a Riemannian manifold with boundary diffeomorphic
to that of the $2$-disc.  Let 
$H(t, \tau): S^1 \times [0,1] \longrightarrow M$ be 
a smooth map such that $H(t,0)=\gamma(t)=\partial M$,
$H(t,1)=p \in M$, and $H(t, \tau)=\gamma_\tau(t)$ is a simple closed curve
parametrized by $t$ for each $\tau \in [0,1]$.  We will say that $H$ is 
a \emph{weakly monotone} homotopy if closed $2$-discs
$D_\tau\subset M$ bounded by $\gamma_\tau$ 
satisfy the inclusion
$D_{\tau_2} \subset D_{\tau_1}$
for every $\tau_1$ and $\tau_2$ with $\tau_1 < \tau_2$.
  If these discs satisfy a stronger condition that  $D_{\tau_2} \subset int D_{\tau_1}$,
where $int D_{\tau_1}$ is the interior of $D_{\tau_1}$, then the
homotopy will be called \emph{strictly monotone}, or just \emph{monotone}.
\end{Def}

Figure \ref*{monotonehomotopy}(a) depicts a strictly monotone homotopy
of $\gamma(t)$ to the point $p$, while
Figure \ref*{monotonehomotopy}(b) depicts a homotopy that is not
monotone. This definition can be trivially extended to homotopies connecting 
a simple closed contractible curve to a point on any closed
Riemannian surface not diffeomorphic to the $2$-disc. There is one technicality, however; 
if $M$ is diffeomorphic to $S^2$, then  there is an ambiguity due to non-uniqueness of $D_{\tau}$.
We agree to resolve it by allowing any possible choice of the system of discs $D_\tau$ that is
continuously dependent on $\tau$, and
that has the monotonicity (or strict monotonicity) property.

\realfig{monotonehomotopy}{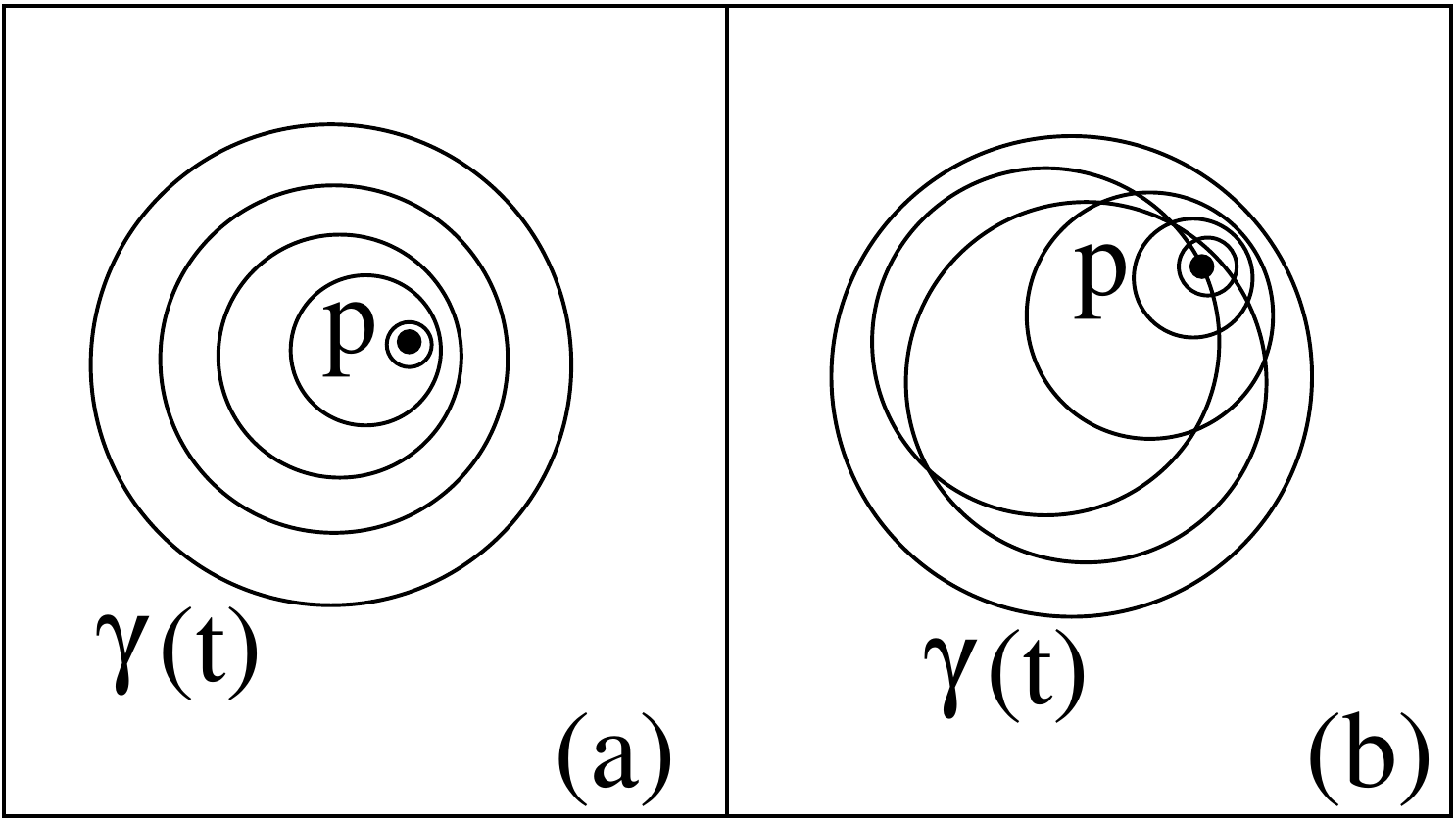}{Monotone and non-monotone homotopies}{0.7\hsize}

There are several natural conjectures concerning the existence of monotone homotopies in the presence
of regular homotopies.  The first one of these conjectures concerns contractions of a disc, and the second
concerns contractions of a simple closed curve on a Riemannian surface.
\begin{Conj}
	\label{conj:hom_to_mono_disc}
	Suppose that $M$ is a Riemannian manifold which is diffeomorphic to a disc.
	If there is a contraction of $\partial M$ through curves of length at most $L$
	then, for every $\epsilon > 0$, there is a monotone contraction of $\partial M$
	through curves of length at most $L + \epsilon$.
\end{Conj}

\begin{Conj}
	\label{conj:hom_to_mono_gen}
	Suppose that $M$ is a Riemannian surface and that $\gamma$ is a simple closed curve
	on $M$.  If $\gamma$ is contractible through curves of length at most $L$
	then, for every $\epsilon > 0$, there is a monotone contraction of a simple closed curve
	$\alpha$ through curves of length at most $L + \epsilon$ which covers $\gamma$.
\end{Conj}

In the second conjecture, when we say that the monotone contraction covers $\gamma$, we mean that $\gamma$
is contained in the the disc with boundary $\alpha$ which is the image of the contraction.

Both of these conjectures can be simplified by using the following theorem of independent interest
proved by the first author and Yevgeny Liokumovich in [CL1]:

\begin{Thm}{(G. R. Chambers, Y. Liokumovich [CL1])} \label{Theoremgy}
Let $M$ be a Riemannian surface.  
Let $\gamma(t): S^1 \longrightarrow M$ 
be a closed curve in $M$.  Suppose there exists a homotopy
$H(t,\tau):S^1 \times [0,1] \longrightarrow M$ such that
$H(t,0)=\gamma(t)$; $H(t,1)=p \in M$, and the length of $\gamma_\tau=H(*,\tau)$ is
at most $L$ for all $\tau \in [0,1]$.  Then, for any $\epsilon >0$, 
there exists a homotopy
$\tilde{H}(t,\tau):S^1 \times [0,1] \longrightarrow M$, such
that $\tilde{H}(t,0)=\gamma(t)$; $\tilde{H}(t,1)=p$; 
the length of $\tilde{\gamma}_\tau=\tilde{H}(*,\tau)$ is at most
$L+\epsilon$, and $\tilde{\gamma}_\tau(t)$ is a
simple closed curve for every $\tau \in [0,1]$.
\end{Thm} 

By applying this theorem, Conjectures \ref*{conj:hom_to_mono_disc} and \ref*{conj:hom_to_mono_gen}
can be reduced to statements about homotopies through simple closed curves.

Our main theorem is a ``gluing" statement about monotone homotopies.  A monotone homotopy $G$
from an initial curve $\gamma_1$ to a final curve $\gamma_2$ traverses an annulus with
$\gamma_1$ as its ``outer" boundary component, and $\gamma_2$ as its ``inner" boundary component.
We begin by stating the relationship between the two annuli formed by the two homotopies
which we are going to glue together:
\begin{Def}
	\label{def:nested}
	Suppose that $G$ and $H$ are two monotone homotopies
	with initial curves $\gamma_1^g$ and $\gamma_1^h$, and with final curves
	$\gamma_2^g$ and $\gamma_2^h$.  Suppose further that two continuous of families of discs
	have been chosen that have boundaries equal to the curves of $G$ and of $H$, respectively.
	For every curve $\gamma$ of $G$ or $H$, we use $D_{\gamma}$ to denote the disc that fills $\gamma$.
	In particular, the curves $\gamma_j^i$ are simple and bound discs $D_{\gamma_i^j}$ for $i \in \{ g, h \}$ and $j \in \{1, 2\}$.
	We will say that $G$ and $H$ are \emph{nested} if $D_{\gamma_2^g} \subset D_{\gamma_1^h}$.
	Furthermore, we will say that $G$ and $H$ are \emph{strictly nested} if
	\begin{enumerate}
		\item	$D_{\gamma_1^h} \subset D_{\gamma_1^g}$
		\item	$D_{\gamma_2^g} \subset D_{\gamma_1^h}$
		\item	$D_{\gamma_2^h} \subset D_{\gamma_2^g}$
	\end{enumerate}
\end{Def}

To state our theorem we also require the following definition:
\begin{Def} \label{simpleintersections}
Let $\beta_1(t): [0,1] \longrightarrow M$ and 
$\beta_2(t):[0,1] \longrightarrow M$ be two closed curves in a 
Riemannian manifold $M$.  

We will say that $\beta_1(t)$ and $\beta_2(t)$
satisfy the \emph{simple intersection property} if, for every two points
of intersection of $\beta_1$ and $\beta_2$, they are consecutive on $\beta_1$
if and only if they are consecutive on $\beta_2$.
\end{Def}

We can now state our main gluing theorem.
\begin{Thm}
	\label{thm:gluing}
	Suppose that $G$ and $H$ are monotone homotopies which are nested, and which
	pass through curves of length at most $L$.  Furthermore, suppose that there is a closed
	curve $\alpha$ such that, using the notation from Definition \ref*{def:nested},
	\begin{enumerate}
		\item	$\alpha$ lies in the closed annulus $cl(D_{\gamma_1^h} \setminus D_{\gamma_2^g})$.
		\item	$\alpha$ minimizes length among all closed curves in this annulus homotopic to $\gamma_1^h$.
		\item	$\alpha$ and $\gamma_1^g$ satisfy the simple intersection property, and $\alpha$ and $\gamma_2^h$
			also satisfy the simple intersection property.
	\end{enumerate}
	Then, for every $\epsilon > 0$, there
	exists a monotone homotopy $K$ with a corresponding family of discs such that
	the disc which fills the initial curve contains $D_{\gamma_1^g}$, and the disc which
	fills the final curve is contained in $D_{\gamma_2^h}$.  Additionally, we
	can construct $K$ so that it is composed of curves of length at most $L + \epsilon$.

	Suppose that $G$ and $H$ are strictly nested, and that both pass through curves of length at most $L$.
	Then, for every $\epsilon > 0$, we can find a monotone homotopy $K$ composed of curves of length
	at most $L + \epsilon$, which begins on $\gamma_1^g$, and which ends on $\gamma_2^h$.
\end{Thm}

This suggests a natural way to approach Conjectures \ref*{conj:hom_to_mono_disc} and
\ref*{conj:hom_to_mono_gen}.  First, we apply Theorem \ref*{Theoremgy} to reduce the problem to the case
where the homotopies pass through simple closed curves.  Next, we choose a finite sequence of simple closed 
curves $\gamma_0, \dots, \gamma_n$ from this homotopy.  We then produce a sequence of 
monotone homotopies $H_0, \dots, H_n$ such that $H_i$ is strictly nested with respect to $H_{i+1}$
for all $i$, and such that $\gamma_i$ is contained in the annulus which $H_i$ traverses.  The idea of
how to produce this sequence is to choose $\gamma_i$ and $\gamma_{i+1}$ close together, and then to perturb each within
its normal neighborhood.  We then apply Theorem \ref*{thm:gluing} to the pair $H_0$ and $H_1$ obtaining a
new homotopy which covers both $\gamma_0$ and $\gamma_1$, and which is strictly nested with respect
to $H_2$.  We then iterate this process for the new homotopy and $H_2$, obtaining a homotopy
which covers $\gamma_0, \gamma_1$, and $\gamma_2$, and which is strictly nested with respect to $H_3$.
Continuing this process for all of the homotopies we achieve the desired result.

This procedure fails because the perturbation argument
does not yield a sequence of \emph{strictly} nested homotopies.  Instead, we can attempt to execute this approach
for nested homotopies.  In this case, the perturbation argument produces a sequence of nested
homotopies, but the proof still fails due to the fact that the inductive step no longer holds.
In particular, if $H_1$, $H_2$, and $H_3$ form a sequence of nested homotopies,
and if we attempted to glue together $H_1$ and $H_2$ as per Theorem \ref*{thm:gluing}, then
the result may no longer be nested with respect to $H_3$.

In \cite{CL2}, the first author and Y. Liokumovich use exactly this technique (and a proof
which mimics the proof of Theorem \ref*{thm:gluing}) to solve conjectures analogous to
Conjecture \ref*{conj:hom_to_mono_disc} and \ref*{conj:hom_to_mono_gen} for one parameter
families of sufficiently regular cycles.  This result plays a crucial role in their proof that every complete manifold of finite
volume contains a (possibly non-compact) minimal surface.

\subsection{Contracting discs in Riemannian manifolds}

Conjectures \ref*{conj:hom_to_mono_disc} and \ref*{conj:hom_to_mono_gen}, if true, have numerous implications.  One which we are particularly interested in
is the following. Assume that $M$ is a smooth surface,
possibly with boundary, and that $\gamma$ is a simple closed curve on $M$. Assume that $\gamma$ can be contracted
to a point via free loops (closed curves) of length at most $L$.  We would like to contract $\gamma$ over closed curves based at a point 
$q \in \gamma$ so that the maximal length of these curves is as small as possible.
Can we estimate the required maximal length in terms of $L$ and the diameter $d$ of $M$?  Such a result would have a large number of applications, some of which
will be discussed at the end of the introduction.  This problem is already interesting when $M$ is a $2$-disc endowed with
a Riemannian metric and $\gamma$ is its boundary.  The simple example of a Riemannian metric
on $M$ that looks like a long thin finger shown in Figure \ref*{monotoneintro1} demonstrates that we cannot estimate the required length in terms of $L$ alone, and this
example suggests that at the very least we need to add a summand equal to $2d$. In this example, we can contract the boundary
of the disc via short closed curves to a point $p$ far from $\partial M$ (see Figure \ref*{monotoneintro1}(a)). To replace this homotopy by one composed of
loops based at a point $q\in\partial M$, we connect $p$ and $q$ by a minimizing geodesic $\tau$. In the course of our new homotopy, we travel along
$\tau$ to one of the closed curves in the original homotopy, travel along this curve, and then return back along $\tau$ (Figure \ref*{monotoneintro1}(b)). At some moment we end up
at a loop that consists of two copies of $\tau$ traversed in opposite directions. This loop can then be contracted to $q$ along itself.

Of course, there are other, more complicated Riemannian metrics
on the $2$-disc, such as the metric depicted in Figure \ref*{monotoneintro2}.  There is also the family of Riemannian metrics considered in [FK].
For these metrics, the connection between 
the length of curves in the ``best'' free loop homotopy and
the length of curves in the ``best'' fixed point homotopy is not so evident.
\par
Our first theorem asserts that if Conjecture \ref*{conj:hom_to_mono_disc} is true,
then adding the summand $2d$ and an arbitrarily small $\epsilon$ to $L$ will always suffice. It is quite possible that one
does not need $\epsilon$, but this does not seem to follow from a compactness argument, as when $\epsilon \longrightarrow 0$, our homotopies can become
wigglier and wigglier. (More formally, we do not establish a bound on the Lipschitz constants of our homotopies as $\epsilon \longrightarrow 0$.)

\begin{Thm} \label{Theoremmain}
Let $M$ be a Riemannian manifold with boundary diffeomorphic
to the standard disc of dimension $2$. Denote its diameter by $d$. Suppose there exists a 
homotopy connecting the boundary $\partial M$ of $M$
to some point $p \in M$  such that
the length of every closed curve in this homotopy does not exceed a real number $L$. 
If Conjecture \ref*{conj:hom_to_mono_disc} is true then,
for any $q \in \partial M$ and 
for any $\epsilon > 0$,
there exists a fixed point homotopy 
that connects $\partial M$ with $q$, 
and passes through loops that are based at $q$ and have length not exceeding $L+2d+\epsilon$.

In particular, if there exists a monotone contraction of $\partial M$ through curves of
length at most $L + \epsilon$, then the result holds.
\end{Thm}
\begin{Pf}{Proof}
By Conjecture \ref*{conj:hom_to_mono_disc}, there exists a strictly monotone homotopy $H$ between
$\gamma(t)$ and $\tilde{p} \in M$ over simple curves of length 
at most $L + \epsilon$.  Fix a point $q$ on $\partial M$, and let $\alpha(s) : [0,1] \rightarrow M$
be a minimal geodesic connecting $q$ to $\tilde{p}$.  The length of $\alpha$ is at most $d$.
For each $\tau \in [0,1]$, there is exactly one $\tau' \in [0,1]$ such that the curve $H(*, \tau')$
goes through $\alpha(\tau)$.  Let this curve be denoted by $H_{\tau}$.  Note that, if $\alpha$ intersects a curve in $H$ multiple times, then
we will be able to find multiple values for $\tau$ that result in the same curve.

Our new contraction $G: S^1 \times [0,1] \rightarrow M$ of $\gamma$ through curves based at $q$ is now defined
as follows.  For each $\tau \in [0,1]$, define $G(*, \tau)$ to be the curve
	$$ \alpha|_{[0, \tau]} * H_{\tau} * \overline{\alpha|_{[0,\tau]}},$$
where $\overline{\alpha|_{[0,\tau]}}$ is the segment of $\alpha$ traversed from $\tau$ to $0$.
Each curve in this homotopy is bounded in length by $L + 2d + \epsilon$.  Furthermore,
it ends at $\alpha|_{[0,1]} * \overline{\alpha|_{[0,1]}}$, which can obviously be contracted
to $q$ through curves based at $q$ of length no more than $2d$.  This completes the proof.
\end{Pf}

Our second theorem deals with the general case of a simple contractible curve on a surface endowed with a Riemannian metric.

\begin{Thm} \label{Theoremmain1}
Let $M$ be a closed Riemannian surface  of diameter $d$.  Let
$\gamma:[0,1] \longrightarrow M$ be a simple
 closed curve in $M$, and $q$ a point on $\gamma$.
 If there exists a homotopy between $\gamma$ and a point that passes through closed curves of length not exceeding $L$, and if Conjecture \ref*{conj:hom_to_mono_gen} is true,
then there exists a homotopy that contracts $\gamma$ to $q$ through loops that are 
based at $q$ and have length $\leq 3L+2d+\epsilon$.

In particular, if there is a monotone contraction of a simple closed curve $\alpha$ through curves of
length at most $L + \epsilon$ which covers $\gamma$, then the result holds.
\end{Thm}

The proof of this theorem is significantly more involved; we will postpone its proof until the end of the article.

\realfig{monotoneintro1}{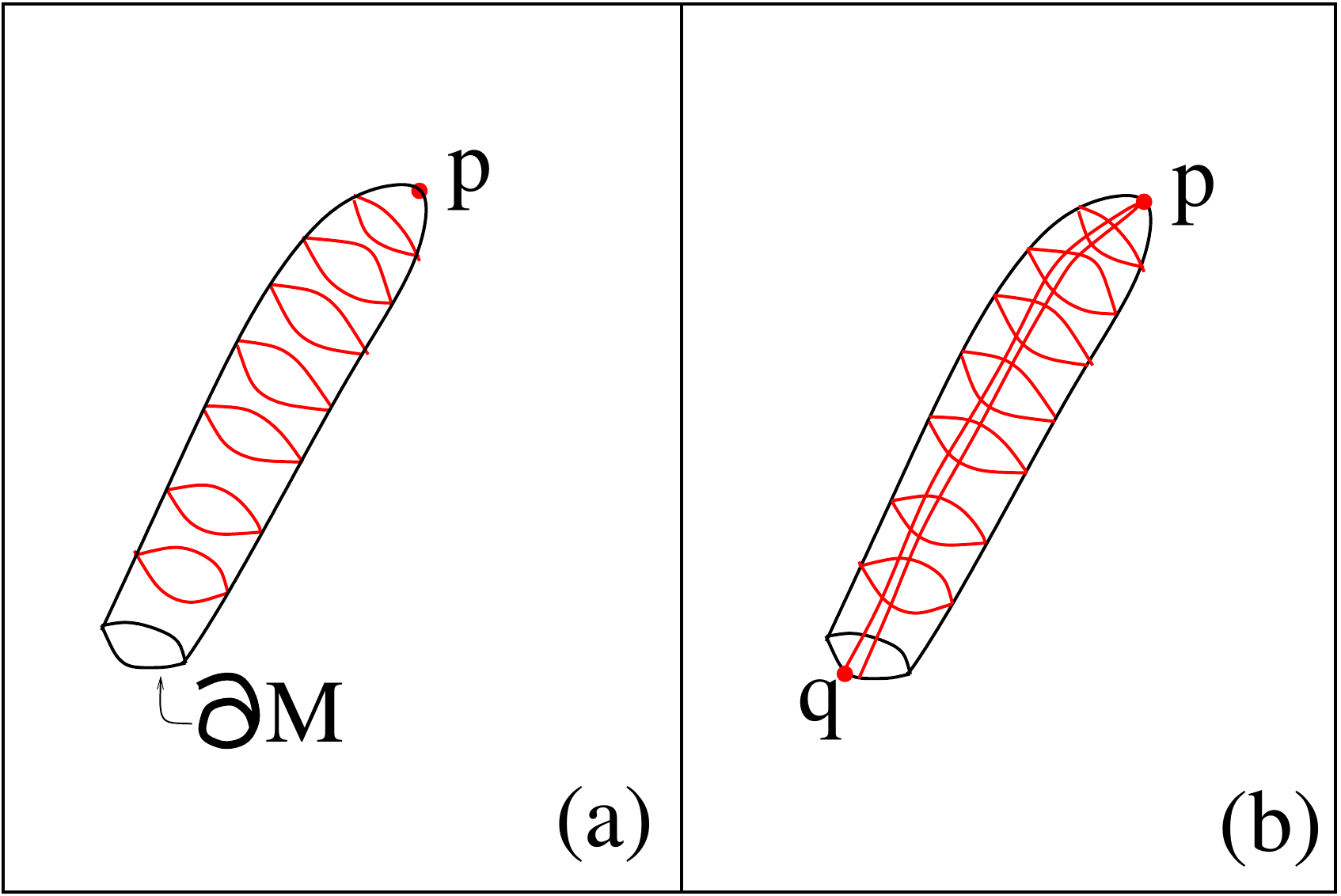}{Long finger}{0.7\hsize}

The questions considered in this paper fall within the realm of questions
of investigating geometric properties of ``optimal'' homotopies.

\realfig{monotoneintro2}{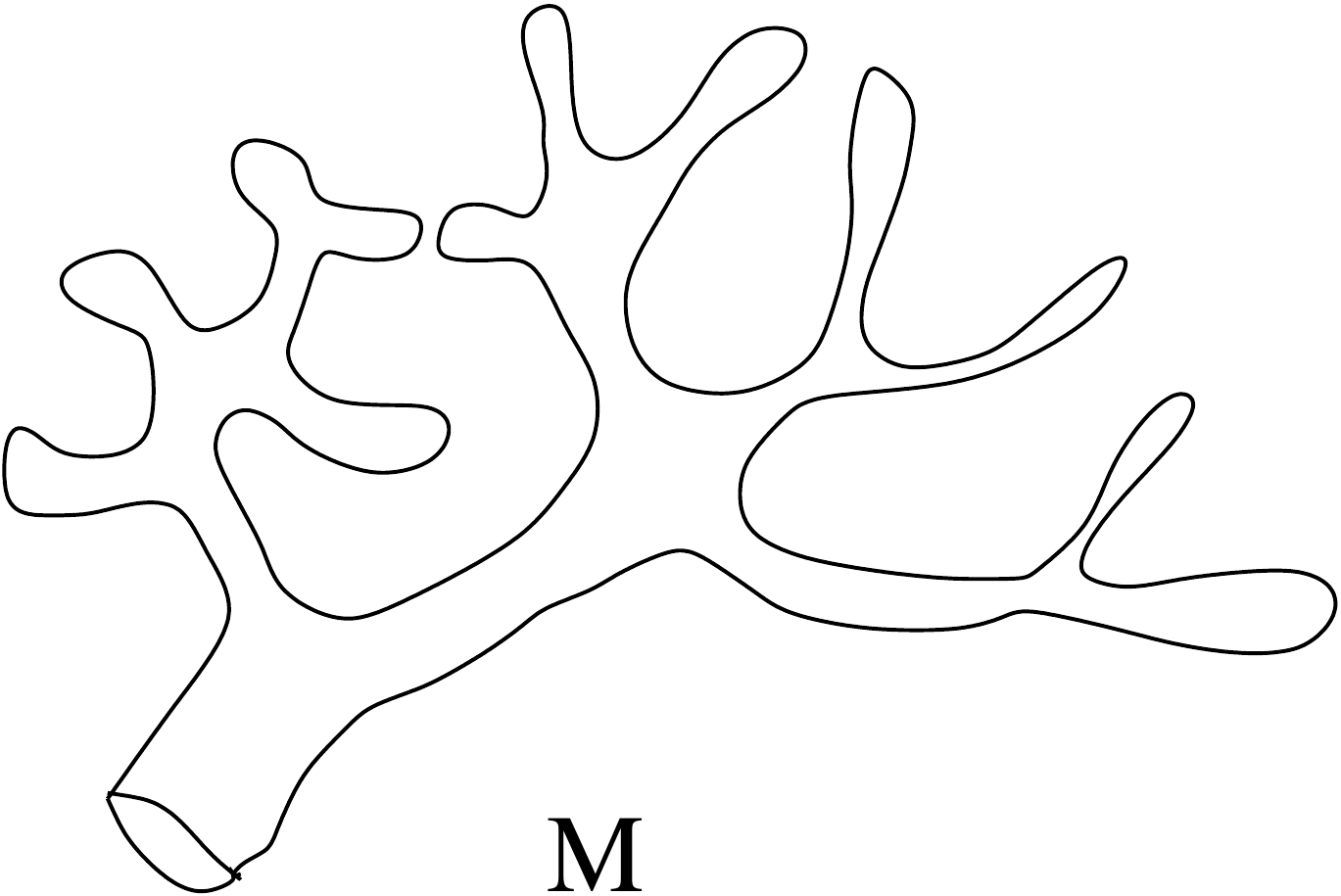}{``Cactus'' metric on the disc}{0.7\hsize}

Another example of such a question was a long-standing question of S. Frankel
and M. Katz posed at the end of their paper [FK].  They asked if one contract the boundary $\partial M$ of a
Riemannian $2$-disc $M$ so that the length of curves in the homotopy
is majorized above in terms of the diameter $d$, area $A$, 
and the length of the boundary of the manifold?
Note that this question is a modification of an earlier question asked by M. Gromov ([Gr], p. 100).
The positive answer to this question was given by 
Y. Liokumovich, A. Nabutovsky and the second author of this paper in [LNR].
In particular, it was shown that $\partial M$ can be contracted over
curves of length at most 
$|\partial M|+200d\max\{1, ln\frac{\sqrt{A}}{d}\}$, where $|\partial M|$ is the length of $\partial M$. This estimate is optimal up to a multiplicative factor in the second term.  When 
$\frac{\sqrt{A}}{d} <<1$, [LNR] provides a better bound of $2|\partial M|+2d+686 \sqrt{A}$.  This was 
improved to the asymptotically tight upper bound $|\partial M|+2d+O(\sqrt{A})$ by P. Papasoglu
in a recent paper [P].
Note that it is impossible to bound the length of curves in 
the best homotopy solely in terms of the area of $M$, as 
an example of a three-legged star  fish with long tentacles depicted
in Figure \ref*{monotoneintro3} demonstrates.  It is also impossible to bound the length
solely in terms of the diameter of $M$, as was proved by 
Frankel and Katz in [FK], answering the original version of the question of M. Gromov
mentioned above.

A closely related family of questions deals with establishing the existence of various upper bounds on
the maximal length of optimal sweep-outs and slicings of surfaces either by closed curves or,
more generally, by cycles. For example, it was shown by Y. Liokumovich that there does not
exist a universal diameter bound for the maximal length of curves or cycles in an
optimal sweep-out of a closed Riemannian surface (see [L1] and [L2]).
On the other hand, F. Balacheff and S. Sabourau have found an upper bound for the maximal
length  of a cycle in an optimal sweep-out of a surface in terms of the genus and the area
of the surface (see [BS]).  Also, a ``short'' sweep-out of a Riemannian $2$-sphere is possible 
if one assumes that there is no ``short'' geodesics of index $0$.  This follows from the results
of C. B. Croke in [C].  In this case, the maximal length of a cycle can be bounded by the area or the diameter of the surface.

\realfig{monotoneintro3}{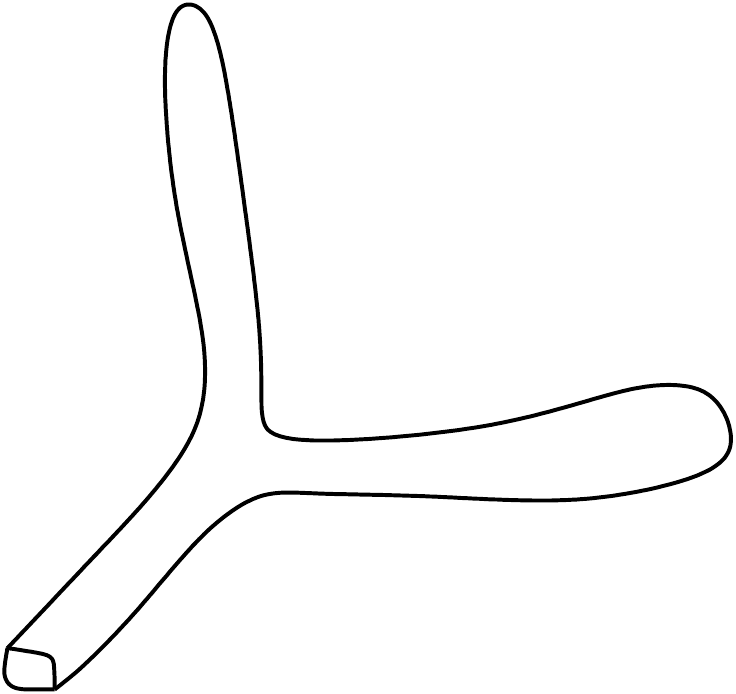}{Three legged star fish}{0.7\hsize}

Note that we cannot hope to prove Conjecture \ref*{conj:hom_to_mono_disc}
when one has a simple curve on a surface, or even a simple closed curve in a disc that is not
assumed to be the boundary of that disc. This fact makes proving Theorem 
\ref*{Theoremmain1} more difficult than Theorem 0.1, and is the reason
for the appearance of the extra $2L$ in our upper bound. We are grateful to Yevgeny Liokumovich
for first attracting our attention to this fact in conjunction with the example shown in
Figure \ref*{notmonotone1}.
This figure depicts a metric on a Riemannian 
$2$-disc and a curve $\alpha_0$ 
such that the optimal homotopy contracting the curve to a point is not
monotone. Notice that there are three bumps depicted in this figure:
two of them are long and thin and the one in the middle is short and asymmetric.
It takes less length to go under the middle bump than over it.
The original curve $\alpha_0$ winds around the two thin bumps,
and goes over the short one.
In order to contract $\alpha_o$ to a point, it has to be stretched
over the thin bumps but, because they are long, the length of the 
curve will necessarily increase in the process.  Thus, to begin with, 
it makes sense to first homotope $\alpha_0$ to $\alpha_1$, which runs below
the middle bump.  $\alpha_1$ is 
shorter than $\alpha_0$, so we can ``spend'' this ``excess'' length 
on dragging the curve over the two thin bumps one at a time.  This
corresponds to the curves $\alpha_2$ and $\alpha_3$ in Figure \ref*{notmonotone1}.
We now have to push the curve over the middle bump.  This homotopy
results in $\alpha_4$, which can then be easily contracted to a point.
The resulting homotopy is not monotone.  Note that $\alpha_0$ is 
not the boundary of the disc, and so this example does not constitute
a counterexample to Conjecture \ref*{conj:hom_to_mono_disc}.

\realfig{notmonotone1}{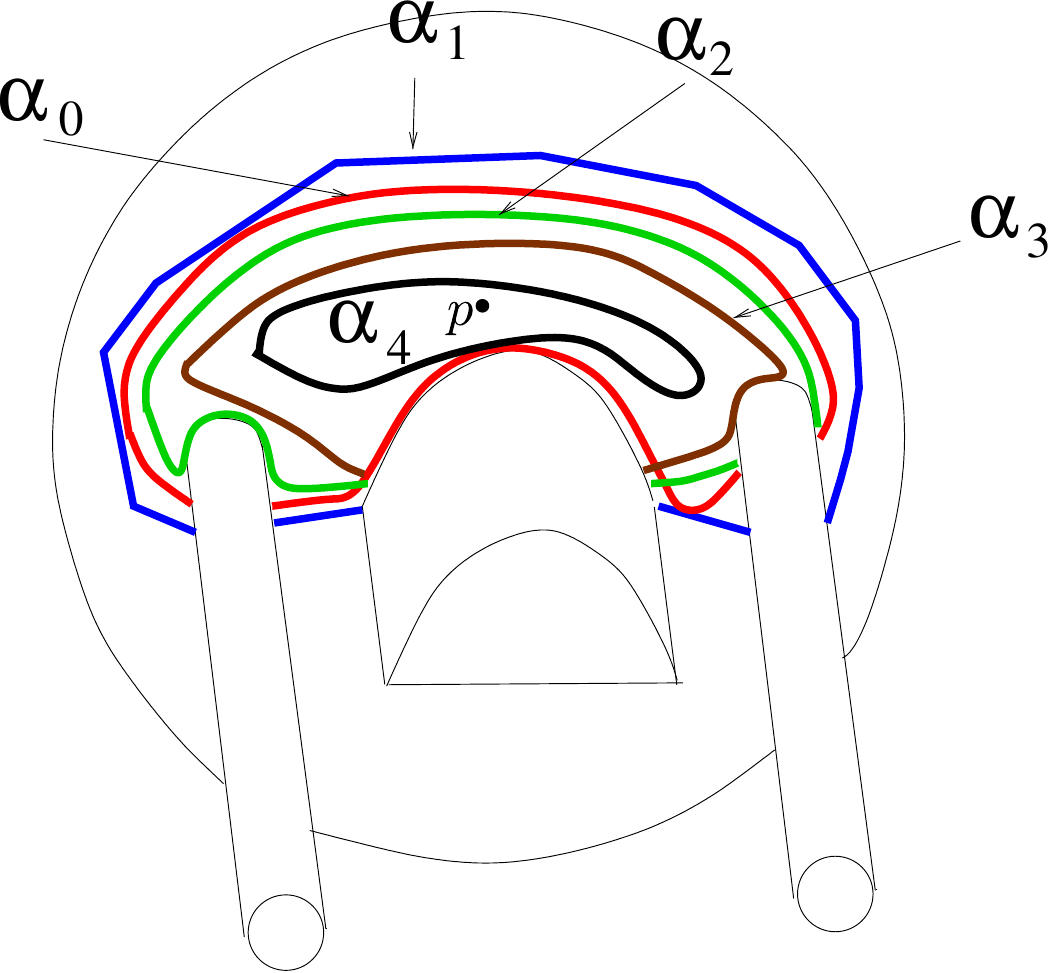}{$\alpha_0$ cannot be contracted to $p$
via a ``short'' monotone homotopy}{0.7\hsize}


\bigskip

\noindent{\bf \large Applications.} If Conjectures \ref*{conj:hom_to_mono_disc} and
\ref*{conj:hom_to_mono_gen} are true, then Theorem \ref*{Theoremmain} and
Theorem \ref*{Theoremmain1} will have 
have numerous immediate and potential applications to the geometry 
of loop spaces of Riemannian $2$-spheres, to questions about
the lengths of geodesics, and to problems about  
optimal sweep-outs.  

In particular, Theorem \ref*{Theoremmain1} provides a canonical way
of obtaining a ``short'' based loop homotopy out of a ``short'' free
loop homotopy on a Riemannian surface.  The second author has encountered
this problem many times, and each time it was solved using ad hoc methods.
Specifically, Theorems \ref*{Theoremmain} and \ref*{Theoremmain1} can be applied 
in the following situations.

\noindent (1) {\bf Lengths of geodesics on Riemannian $2$-spheres.}

Let $p,q$ be an arbitrary pair of points on a Riemannian $2$-sphere
$M$ of diameter $d$.  A. Nabutovsky together with the second author
have demonstrated that there exist at least $k$ geodesics joining them
of length at most $22kd$ (see [NR2]). If $p = q$, then this bound becomes 
$20kd$ (see [NR1]).  We have noticed that applying 
Theorem \ref*{Theoremmain1} dramatically decreases the complexity 
of proofs in [NR1] and [NR2], and improves the bounds in [NR2] to $16kd$, and
the bounds in [NR1] to $14kd$. 
These improvements are due to the fact that the main technical
difficulty encountered in [NR1] and especially in [NR2] is the possible formation of
intersections between various closed curves in a homotopy between a closed curve and a point.

\noindent (2) {\bf Geometry of the loop spaces of Riemannian $2$-spheres.}
Applying Theorem \ref*{Theoremmain1} will immediately 
generalize Theorem 1.1 in [NR3] to the {\it free} loop space of
a Riemannian $2$-sphere $M$.  To be more precise, one can show that
any map $f:S^m \longrightarrow \Lambda M$, where $\Lambda M$ is a free loop
space on $M$, is homotopic to a map $\tilde{f}:S^m \longrightarrow \Lambda M$
that passes through curves of length that do not exceed $L=L(m,k,d)$.  Here, 
$d$ is the diameter of $M$, $k$ is the number of distinct 
non-trivial periodic geodesics
on $M$ of length at most $2d$, and $L$ is a function of $m, k, d$ that
can be written down explicitly. Moreover, one can explicitly majorize the 
lengths of loops in an ``optimal'' homotopy connecting $f$ and $\tilde{f}$
in terms of $n,k,d$ and $\sup_{x \in S^m} length(f(x))$.

\section{Proof of Theorem \ref*{thm:gluing}}

In this section we prove Theorem \ref*{thm:gluing}, which we recall below:
\begin{Thm*}

	Suppose that $G$ and $H$ are monotone homotopies which are nested, and which
	pass through curves of length at most $L$.  Furthermore, suppose that there is a closed
	curve $\alpha$ such that, using the notation from Definition \ref*{def:nested},
	\begin{enumerate}
		\item	$\alpha$ lies in the closed annulus $cl(D_{\gamma_1^h} \setminus D_{\gamma_2^g})$.
		\item	$\alpha$ minimizes length among all closed curves in this annulus homotopic to $\gamma_1^h$.
		\item	$\alpha$ and $\gamma_1^g$ satisfy the simple intersection property, and $\alpha$ and $\gamma_2^h$
			also satisfy the simple intersection property.
	\end{enumerate}
	Then, for every $\epsilon > 0$, there
	exists a monotone homotopy $K$ with a corresponding family of discs such that
	the disc which fills the initial curve contains $D_{\gamma_1^g}$, and the disc which
	fills the final curve is contained in $D_{\gamma_2^h}$.  Additionally, we
	can construct $K$ so that it is composed of curves of length at most $L + \epsilon$.

	Suppose that $G$ and $H$ are strictly nested, and that both pass through curves of length at most $L$.
	Then, for every $\epsilon > 0$, we can find a monotone homotopy $K$ composed of curves of length
	at most $L + \epsilon$, which begins on $\gamma_1^g$, and which ends on $\gamma_2^h$.
\end{Thm*}

\begin{Pf}{Proof}
Without loss of generality, let us assume that
the homotopies $G$ and $H$ are strictly monotone
and pass through smooth curves.  Let $D_g$ be the disc corresponding to the final curve of
$G$, and let $D_h$ be the disc corresponding to the initial curve of $H$.  We may also assume that $\partial D_g$ and $\partial D_h$ are disjoint.
If necessary, all of these properties can be  achieved
by an arbitrarily small perturbation of $G$ and $H$.

Let us consider the closure of the annulus $D_h \setminus D_g$.
Let $\alpha$ be the shortest closed 
curve  
among all closed curves in this closed annulus that are homotopic to $\partial D_h$.  By hypothesis,
we may assume that $\alpha$ and the final curve of $G$ have the simple intersection property,
and that $\alpha$ and the initial curve of $H$ also have the simple intersection property.
Note that
the length of $\alpha$ is at most $L$.
It is easy to see that $\alpha$ is a simple closed curve.
Let $D_\alpha$ be a closed domain
diffeomorphic to the $2$-disc that has $\alpha$ as its boundary and is contained
in $D_h$.  We may assume that the initial and final curves in $G$ and $H$ are transverse to $\alpha$ up to a small
perturbation.

We will form $K$ as follows.  We will define a process which transforms $G$
into a weakly monotone homotopy which starts at the boundary of a disc which contains $D_{\gamma_1^g}$,
and which ends at $\alpha$, which transforms $H$ into a weakly monotone homotopy which starts
at $\alpha$ and ends on the boundary of a disc contained in $D_{\gamma_2^h}$.  This process will increase the lengths of curves by at most $\epsilon$.
$K$ will then will then be defined as the concatenation of the
modified $G$ and $H$ perturbed slightly so that it is strictly monotone.

We will define this process in three steps.  Step 1 defines the above process for modifying
$G$, Step 2 defines the above process for modifying $H$, and Step 3 defines the concatenation of the
two homotopies.

In all of these steps, the homotopies formed may not be
strictly monotone, and may contain curves with tangential self-intersections (or segments which
agree).  These will be resolved in the final step by a slight perturbation.  Throughout
this section, we will use the notation $D_{\gamma, x}$ to denote the disc that has boundary
equal to the curve $G(t,x)$, and we will use the notation $D_{\beta,x}$ to denote the disc that
has boundary $H(t,x)$.

\noindent{\bf Step 1.} We will modify the homotopy $G(t,x)$ to
obtain a new homotopy $\tilde{G}(t,x)$ so that the new
homotopy will be monotone,  $\tilde{G}(t,0)$ will be outside of the disc bounded by $G(t,0)$, 
and $\tilde{G}(t,1)=\alpha(t)$.

First, we will construct a $1$-parameter family of simple closed curves $\tilde{\gamma}_x(t)$
that can possibly have a finite number of discontinuities (as a function of the parameter $x$). 
We will then remove the discontinuities.

The new curves $\tilde{\gamma}_x(t)$ will be constructed
as follows.  
The general principle is to ``push" (arcs of) $\gamma_x$ outside of the interior
of the disc $D_\alpha$.
More specifically, if $D_{\gamma,x} \subset D_\alpha$ as in 
Figure \ref*{monotoneg1}(a),
then we will let $\tilde{\gamma}_x(t)=\alpha(t)$.
If $D_\alpha \subset D_{\gamma,x}$ or $D_\alpha \cap D_{\gamma,x} =\emptyset$ as
in Figure \ref*{monotoneg1}(b), 
then we will let $\tilde{\gamma}_x(t)=\gamma_x(t)$.

\realfig{monotoneg1}{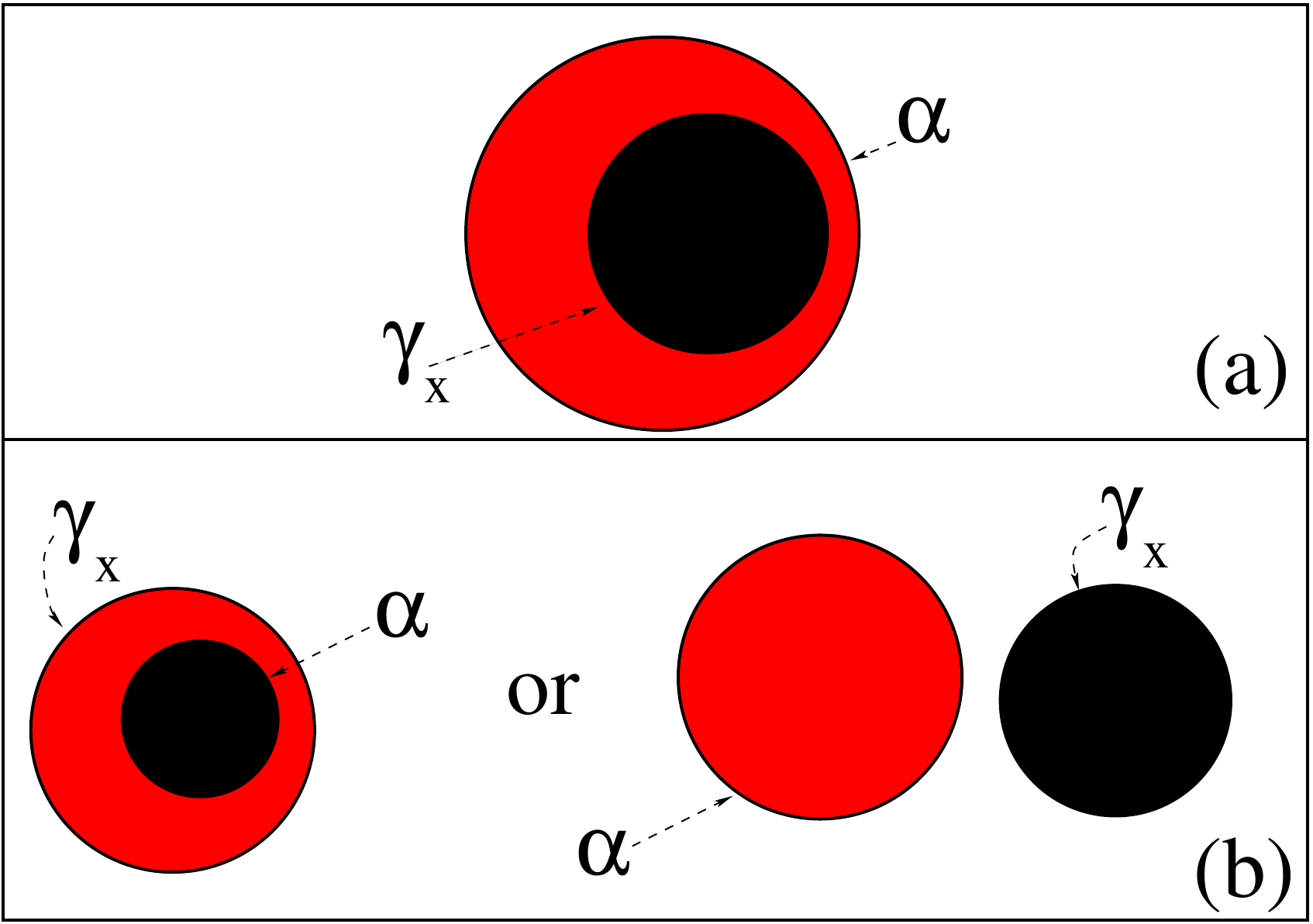}{$\gamma_x$ relative to $\alpha$}{0.7\hsize}

Now suppose that $D_{\gamma,x} \cap D_\alpha \neq \emptyset$ and that neither one
is a subset of the other.  Let us consider $0=t_0<t_1<...<t_n=1$, 
a subdivision of $[0,1]$ such that $\gamma_x(t_i)=\alpha(s_{j(i)})$ and
$\gamma_x \neq \alpha$ otherwise.

Let us consider those arcs of $\gamma_x$ that are inside
$D_\alpha$. Clearly those arcs will be inside $D_h$ as well.
However, also, by the monotonicity of the homotopy $G$, these arcs will be outside
of $D_g$. Thus, they will lie in 
the closed  annulus $cl(D_h-D_g)$, the closure of $D_h-D_g$.

Consider one such arc of $\gamma_x$ between $\gamma_x(t_i)$ and $\gamma_x(t_{i+1})$.
These two points coincide with points $\alpha(s_{j(i)})$ and $\alpha(s_{j(i+1)})$ which subdivide
the curve $\alpha$ into two arcs.  Let us select the arc $A$ of $\alpha$
between these two points which
is path homotopic to the arc of $\gamma_x$ with the same endpoints which lies
inside the closed annulus between $D_h$ and $D_g$.
Now replace this arc of $\gamma_x$ by $A$. The resulting
curve will not be contractible in the closed annulus between $D_g$ and $D_h$.
Also, the length of $A$ does not exceed the length of the
arc of $\gamma_x$ that it replaced, as otherwise $\alpha$ would not be the shortest non-contractible curve in the closed annulus.
Therefore, the resulting
curve will have length bounded by the length of 
the original curve. 
Note that $\tilde{G}(t,0)$ is outside of $G(t,0)$ since $\gamma^g_2$ and $\alpha$ satisfy
the simple intersection property, and so 
(arcs of) curves in the homotopy $G$ can move only outside. Furthermore,
if $G$ and $H$ are strictly nested, then $\gamma^g_2$ and $\alpha$ do not
intersect, and so $\tilde{G}(t,0) = G(t,0)$.
Also, by the hypotheses of lemma,
$D_{\gamma,1} \subset D_\alpha$. As a result, $\tilde{G}(t,1)=\alpha$.

Additionally, note that the resulting map $\tilde{G}(t,x)$ regarded as a $1$-parametric
family of closed curves depending on the parameter $x$ is monotone,
but its dependence on $x$ is not necessarily continuous.  

The possible discontinuities may only occur at those
curves $\gamma_x$ that are tangent to $\alpha$ at some points.
(Here it is convenient for us to think that $\alpha$ is smooth. While, in 
general, $\alpha$ will be only piecewise smooth, we can perturb it into
an arbitrarily close smooth curve of arbitrarily close length. This replacement
of $\alpha$ by a very close smooth curve adds a summand to our estimates 
that can be made 
arbitrarily small. Therefore, without loss of generality
we can assume that $\alpha$ is smooth.)
Without loss of generality, we can also assume that, for each value of $x$,
there is at most one point where $\gamma_x$ and $\alpha$ are tangent,
and that all of these tangencies are non-degenerate.
(If not,
we can perturb $G$ so that the resulting homotopy has these properties, and so that the curves increase in length only by an arbitrarily small amount.) 
Let us modify the family $\tilde{G}(t,x)$ so that the singularities at tangential points are resolved.

In order to do this, let us consider $x_0$ such that
the curve $\gamma_{x_0}$ is
tangent to $\alpha$ at the point $\gamma_{x_0}(s_{x_0})$.  Let us consider the following
two cases:

\noindent {\bf Case A.}  The set of intersections of 
$\gamma_{x_0}$ with the curve $\alpha$ consists of the one point 
$\gamma_{x_0}(t_{x_0})=\alpha(s_{x_0})$. 
First, consider the situation in which $int D_{\gamma,{x_0}} \cap int D_\alpha
=\emptyset$ as in 
Figure \ref*{monotonecase1}(a),  where $int$ denotes the interior of the set.  This situation is not possible since the homotopy $G$ is monotone, and since $D_g \subset D_\alpha$.
Hence, there are two possibilities: the
curve $D_{\gamma,{x_0}} \subset D_\alpha$, as in 
Figure \ref*{monotonecase1}(b), or $D_\alpha \subset D_{\gamma,x_0}$, as 
in Figure \ref*{monotonecase1}(c).
In the former case, we let $\tilde{\gamma}_{x_0}=\alpha$, and
in the latter case, we let $\tilde{\gamma}_{x_0}=\gamma_x$.

\realfig{monotonecase1}{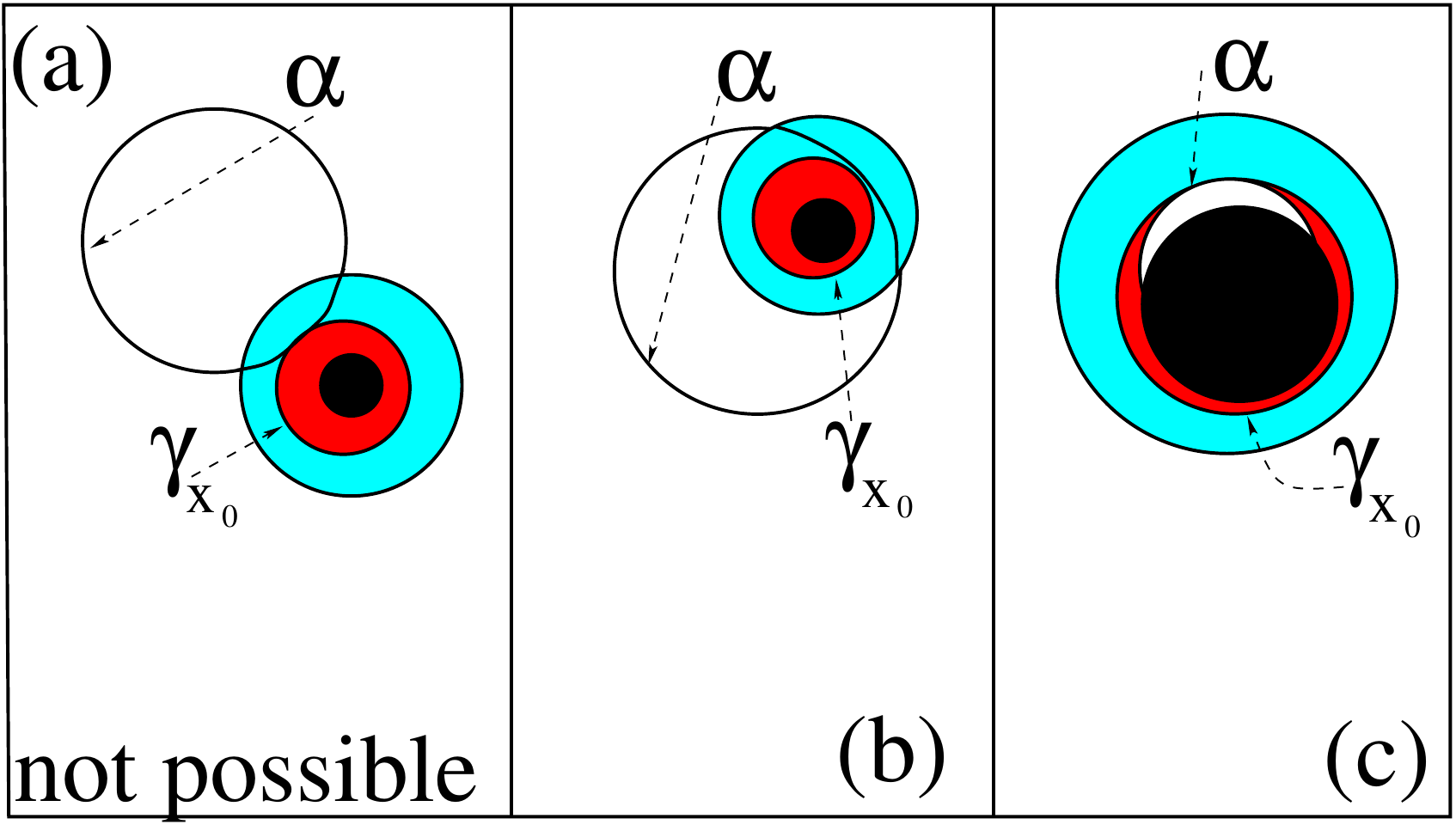}{One intersection between $\alpha$ and $\gamma_{x_0}$}{0.7\hsize}

In both of these cases, the problem with continuity
of the newly constructed $1$-parametric family of curves
at $x_0$ cannot arise.  
Indeed, let us consider the first case, when $D_{\gamma, {x_0}} \subset D_{\alpha}$ as in Figures \ref*{monotonecase1}(a)
and \ref*{monotonecase1a}(b).  The second situation will be 
analogous.  Let us consider a family of curves 
$\gamma_x(t), x \in (x_0-\delta, x_0+\delta)$.  We can choose 
$\delta$ to be small enough so that there are no 
points where  $\gamma_x$ and $\alpha$ are tangent for any $x\not= x_0$ in
the interval $(x_0-\delta,x_0+\delta)$. 
By the monotonicity of the homotopy $G$,
$D_{\gamma,x} \subset D_\alpha$ if $x \in (x_0,x_0+\delta)$, so 
$\tilde{\gamma}_x=\alpha$ for each  $x \in (x_0+\delta)$, as depicted in Figures \ref*{monotonecase1a}(b) and \ref*{monotonecase1a}(c).  
Since $\gamma_{x_0} = \alpha$, $\tilde{\gamma_x}$ is continuous at $x_0$ from the right.  For the other direction, each curve
$\gamma_x(t)$ with $x \in (x_0-\delta, x_0)$ has exactly two points of intersection
with $\alpha$ 
that we will denote $\gamma_x(t_1^x)$ and $\gamma_x(t_2^x)$, where $t_1^x<t_2^x$. 
Note that $\gamma_x(t_1^x)=\alpha(s_1^x)$ and $\gamma_x(t_2^x)=\alpha(s_2^x)$, 
where $s_1^x<s_{x_0}<s_2^x$ by the monotonicity of the homotopy $G$, and both 
$\alpha|_{[s_2^x,s_1^x]}$ and $\gamma_x|_{[t_2^x,t_1^x]}$ 
vary continuously with $x$
as $x$ approaches $x_0$ from the left, and
approach $\alpha(s_{x_0})$ and $\gamma_{x_0}(t_{x_0})$, respectively.
Note that the segment 
 $\gamma_x|_{[t_2^x,t_1^x]}$ is inside of $D_\alpha$.  Hence, our algorithm continuously
replaces this segment by the segment $\alpha|_{[s_2^x,s_1^x]}$ .
Thus, we obtain continuity from the left, as shown in Figure \ref*{monotonecase1a}(d).

\realfig{monotonecase1a}{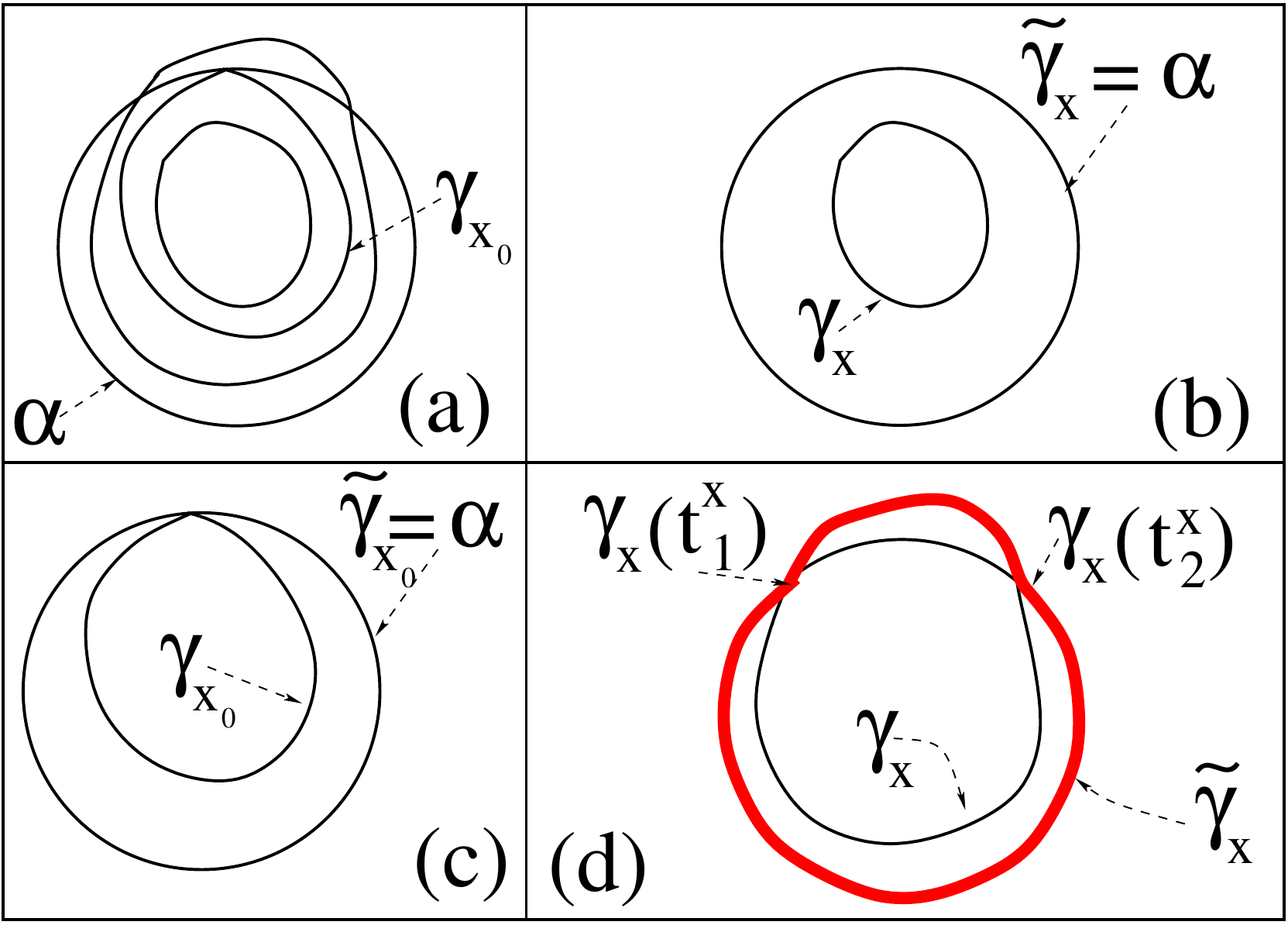}{$\tilde{\gamma}_x$ for $x$ near $x_0$}{0.7\hsize}

\noindent {\bf Case B.}
Let us now assume that $\gamma_{x_0}(t)$ intersects $\alpha$ at 
$2k+1$ points for some $k \geq 1$.  Let $0 \leq t_1<...<t_{2k+1}\leq 1$ be
a partition of the unit interval such that $\gamma_{x_0}(t_j)=\alpha(s(t_j))=
\alpha(s_j)$,
and where all the intersections are transverse except for 
$\gamma_{x_0}(t_{j_0}))$, where it touches the curve $\alpha$ at
$\alpha(s_{j_0})$.

Let us consider the arcs $a=\gamma_{x_0}|_{[t_{j_0-1},t_{j_0}]}$ and 
$b=\gamma_{x_0}|_{[t_{j_0}, t_{j_0+1}]}$.  There are two possibilities to consider:
either both $a$ and $b$ lie inside $D_\alpha$, or both of them lie
outside.  It is not possible for one of the arcs to be inside,
and for the other to be outside, because that would imply that 
either the intersection at $\gamma_{x_0}(t_{j_0})$ is transverse, or
the tangency at this point is degenerate, and either of these options contradicts
our assumptions.

First, assume that both arcs $a$ and $b$ are outside of $D_\alpha$.  Then, according to our algorithm, 
the arcs of the new curve $\tilde{\gamma}_{x_0}$ on the
interval $[t_{j_0-1}, t_{j_0+1}]$ will remain unchanged.  Let us
consider $\gamma_x(t), x \in (x_0-\delta, x_0+\delta)$, where 
$\delta$ is selected so that the curves in the family do not have
any additional tangential intersections, as in Figures \ref*{monotonecase2a}(a) and \ref*{monotonecase2a}(b).
Then, for any nearby curve, the corresponding arc
$\gamma_x|_{[t^x_{j_0-1}, t^x_{j_0+1}]}$ is either outside $D_\alpha$, 
or has two more additional transverse
intersections with $\alpha$ in small neighbourhood of 
$\gamma_{x_0}(t_{j_0})$.
Let us denote these intersections as $t^x_*<t_{j_0}$ and $t^x_{**}>t_{j_0}$.
Thus, $\tilde{\gamma}_x$ will locally be formed by continuously replacing $\gamma_x|_{[t^x_{*}, t^x_{**}]}$ by an arc
of $\alpha$ between the same endpoints.  Hence, in this case, 
$\tilde{\gamma}_x$ changes continuously when $x$ is near $x_0$.  
Note that while $s_{j_0-1}, s_{j_0},$ and $s_{j_0+1}$ can come in any order
as is indicated in Figure \ref*{monotonecase2a}, this fact does not affect
the above analysis.

\realfig{monotonecase2a}{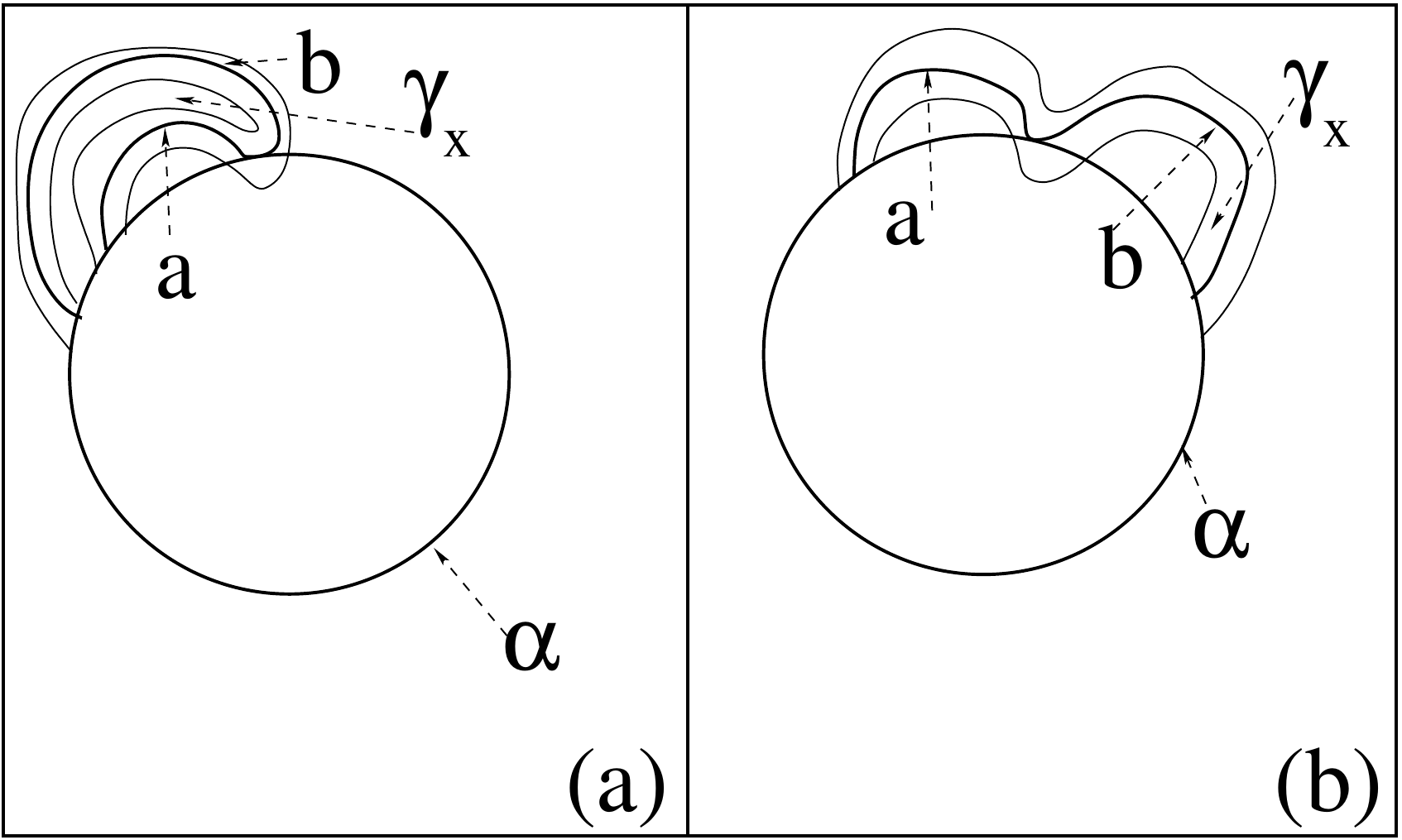}{Outside arcs}{0.7\hsize}

Thus, the only problematic case is when the arcs $a$ and $b$ are inside
$D_\alpha$ as in Figures \ref*{monotonecase2b} and \ref*{monotonecase2c}.  
Let $m,l,$ and $k$ be three points of intersection between 
$\gamma_{x_0}$ and $\alpha$, where $l$ is the tangential point, and 
$m$ and $k$ are its neighbours. By ``neighbours'' we will mean points of intersection
that are the closest to the tangential point along
$\gamma_{x_0}$.  That is, if $l=\gamma_{x_0}(t_{j_0})$, 
then $m=\gamma_{x_0}(t_{j_0-1})$ and $k=\gamma_{x_0}(t_{j_0+1})$.
It is, however, quite possible that along $\alpha$ there are intersection
points that are closer to $l$ than $m$ and $k$ as 
in Figure \ref*{monotoneintersection}.

\realfig{monotoneintersection}{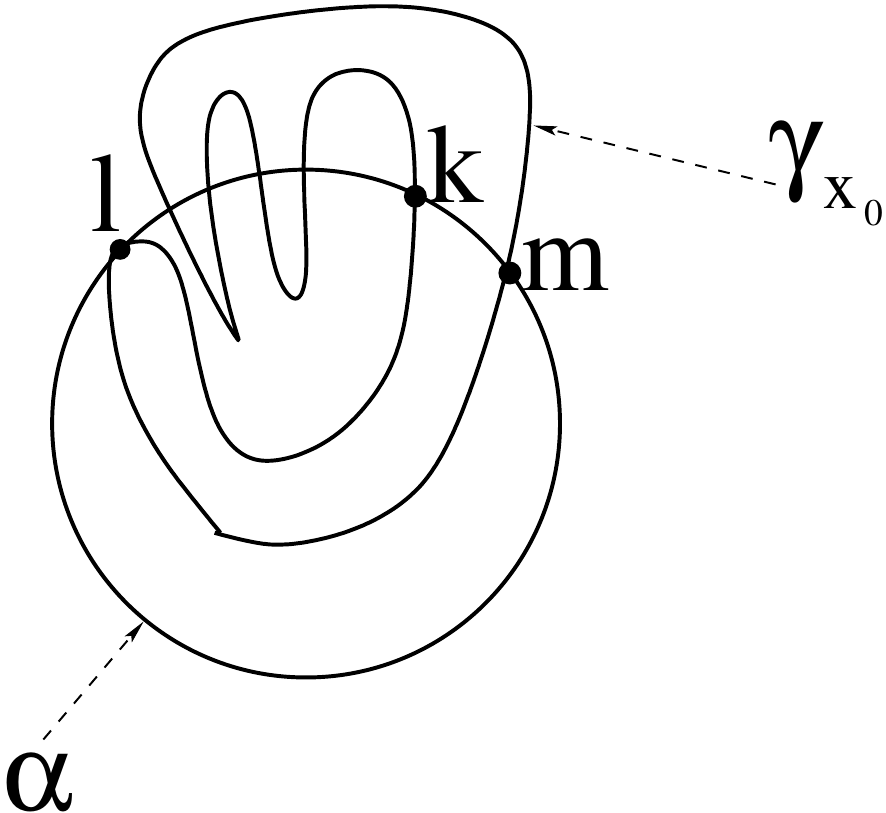}{Complicated intersections between $\alpha$ and $\gamma_{x_0}$}{0.7\hsize}

Each pair of points subdivides the curve 
$\alpha$ into two segments connecting them.
Let us denote the two arcs connecting $m$ with $l$ as $a_{ml}$ and
$\tilde{a}_{ml}$, the two arcs connecting $m$ and $k$ as $a_{mk}$ and
$\tilde{a}_{mk}$, and the two arcs connecting $l$ and $k$ as $a_{lk}$ and
$\tilde{a}_{lk}$.  Also, let us select $a_{ml}, a_{mk}$, and $a_{lk}$ so
that their interiors do not intersect.  Without loss of
generality, let $a$ be the arc of $\gamma_{x_0}$ between
$m$ and $l$ in $D_\alpha$, while $b$ is the arc between $l$ and $k$ in $D_\alpha$.

Since $a$ and $b$ are inside $D_\alpha$, our algorithm requires
that we change them to the corresponding arcs of $\alpha$.
The possible discontinuity is a result of the following situation.  Let us define
two options for replacing arcs of $\gamma_{x_0}$ with arcs of $\alpha$.
Note that, as before, the arc of $\alpha$ that we are replacing an arc of $\gamma_{x_0}$ with is the one that is path homotopic to the arc of $\gamma_{x_0}$ inside the closed
annulus $cl(D_h \setminus D_g)$.

\noindent {\bf Option 1.} Separately replace $a$ and $b$ with the appropriate arcs of $\alpha$.

\noindent {\bf Option 2.} Let $c=a*b$.  Replace $c$ with the appropriate arc of $\alpha$.

As we will see, the two different options will some time result
in the same curve, and some times not.  Our algorithm, at point $x_0$, always uses
Option 1.  Let us, however, consider the curves that are formed by both options.

Let us first consider Option 1.  Without loss of generality, 
suppose that arc $a$ is changed to $a_{ml}$.  (The other option
is $\tilde{a}_{ml}$.) Now there are 
two possibilities for the arc $b$.  It will either be changed to 
$a_{lk}$ or to $\tilde{a}_{lk}$.

Now let us consider Option 2.  Note that if $b$ was
changed to $a_{lk}$, then 
$c$ must be changed to $\tilde{a}_{mk}=a_{ml}*a_{lk}$.
Thus, in this case, it does not matter whether we used Option 1 or 
Option 2 (see Figure \ref*{monotonecase2b}).

However, if $b$ was changed to $\tilde{a}_{lk}$, then 
$c$ must be changed to $a_{mk}$.  
While $a_{mk} \neq a_{ml}*\tilde{a}_{lk}$, we have that 
$a_{ml}*\tilde{a}_{lk}=a_{ml}*\bar{a}_{ml}*a_{mk}$.
(Recall that $\bar{a}_{ml}$ denotes $a_{ml}$ traversed in the opposite direction,
from $l$ to $m$.)
Thus, $a_{mk}$ is path homotopic to $a_{ml}*\tilde{a}_{lk}$ by 
simply contracting $a_{ml}*\bar{a}_{ml}$ to $m$ along itself.
Observe that the length of the curve during this homotopy changes monotonically (see Figure \ref*{monotonecase2c}).

One can see that while the former situation does not create
a discontinuity, the latter situation does.  
Let $\delta>0$ be once again small enough so that
$\gamma_x$ and $\alpha$ do not have any additional tangential 
points on the interval $(x_0-\delta, x_0+\delta)$.

It is possible that, when $x' \in (x_0-\delta, x_0)$, 
$\tilde{\gamma}_{x'}$ will approach the curve obtained using
Option 2 as $x'$ approaches $x_0$, and for $x'' \in (x_0, x_0+\delta)$,
$\tilde{\gamma}_{x^{''}}$ will approach the curve obtained using Option 1
as $x^{''}$ approaches $x_0$ (or the other way around).  This is shown in Figure \ref*{monotonecase2d}.

\realfig{monotonecase2b}{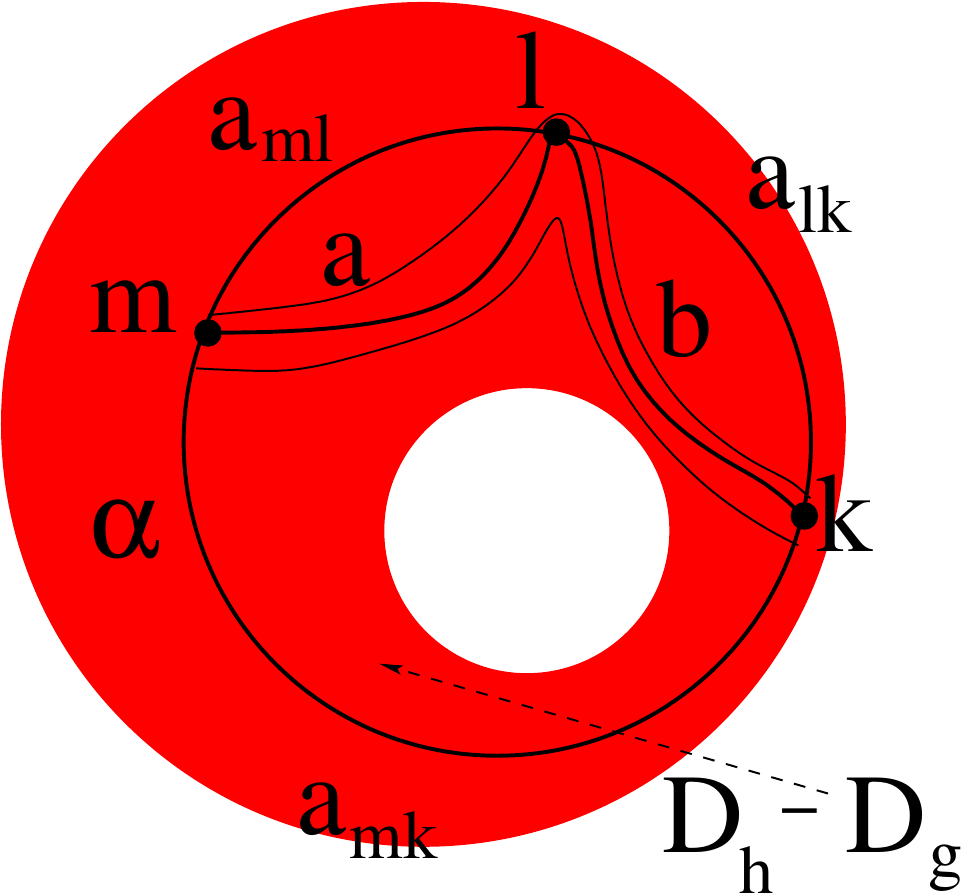}{Inside arcs}{0.7\hsize}

\realfig{monotonecase2c}{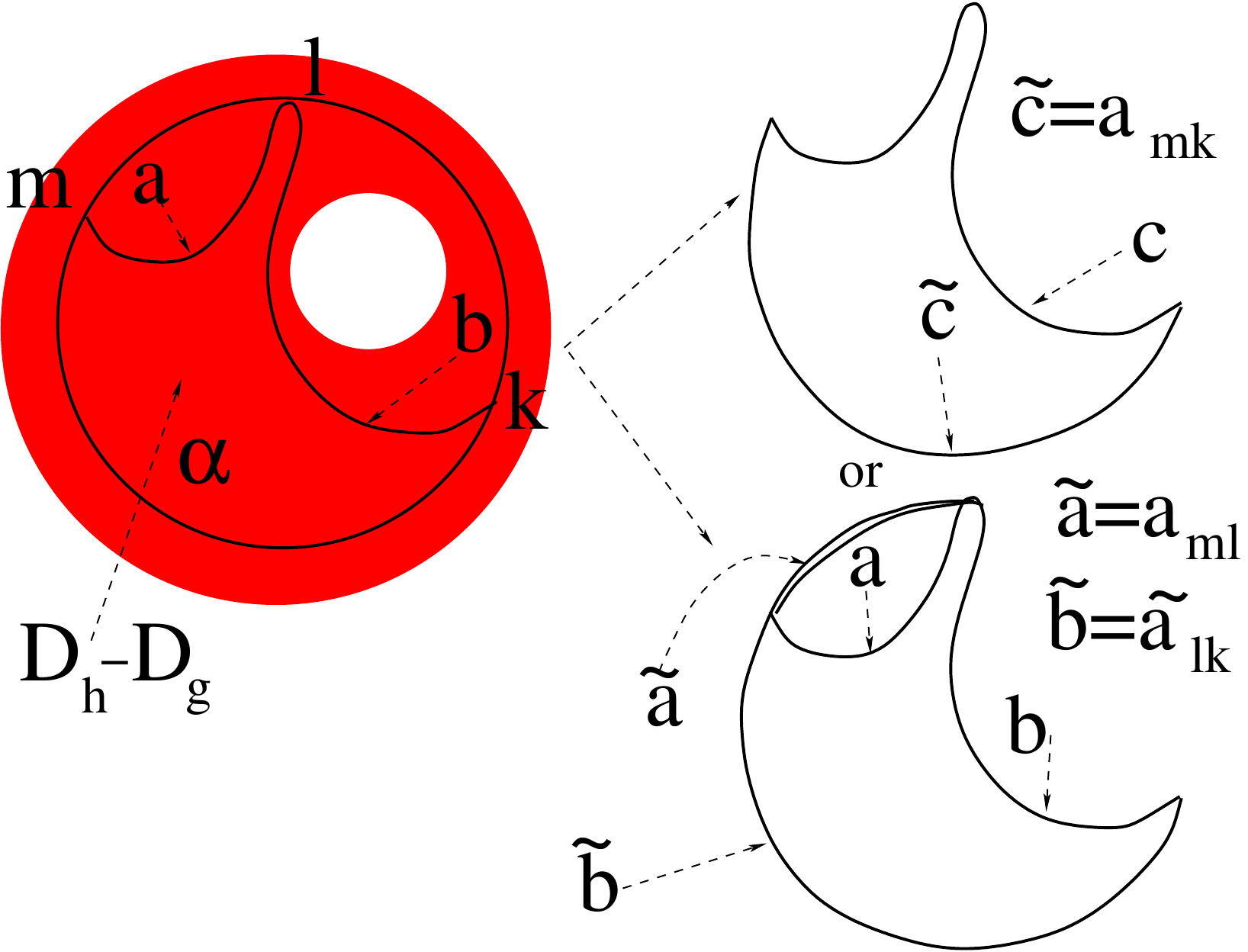}{Two ways to change
arcs}{0.7\hsize}

\realfig{monotonecase2d}{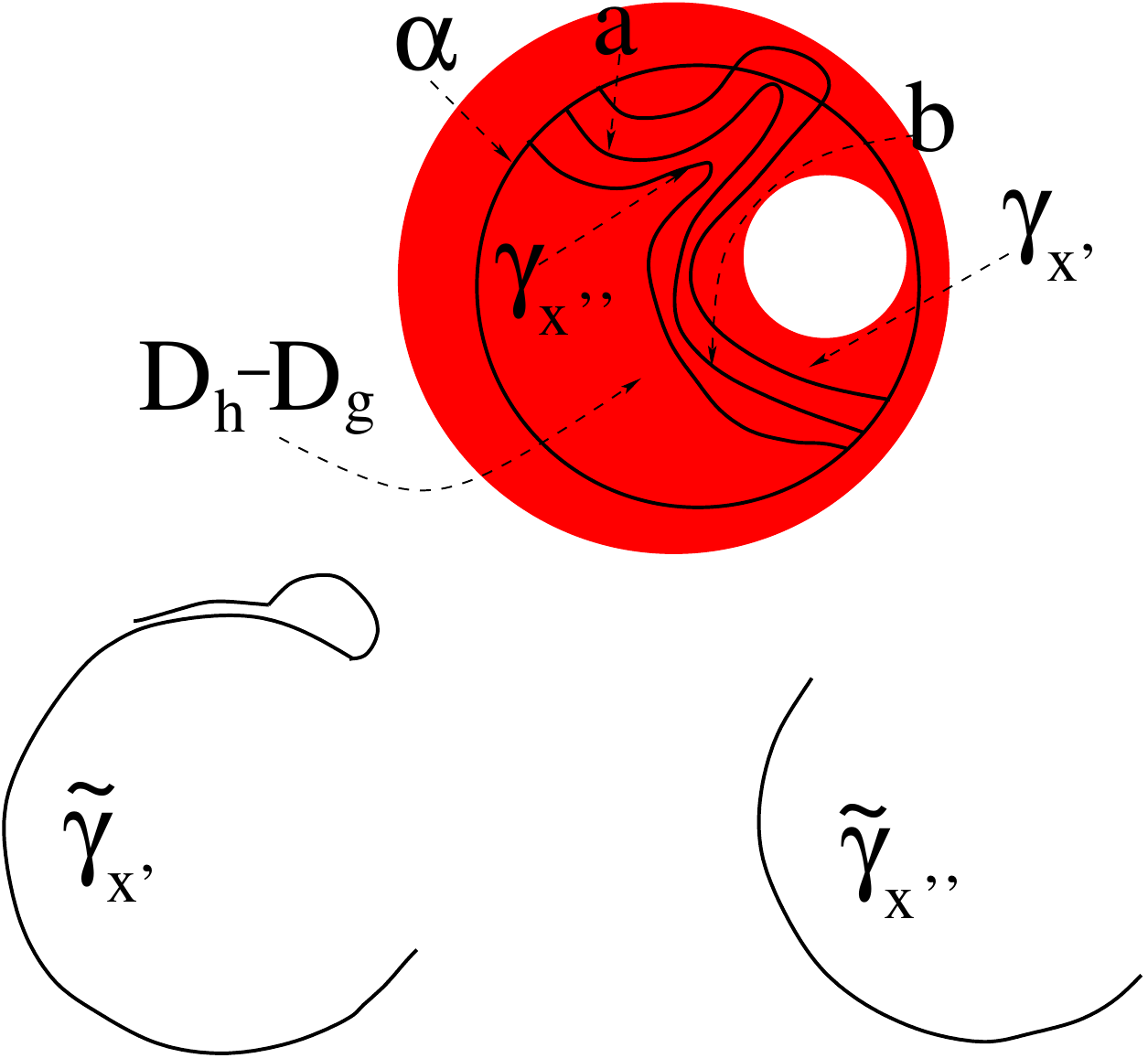}{Homotopy between the curves formed by Option 1 and Option 2}{0.7\hsize}

If we include the homotopy between the curves formed by Option 1 and Option 2 as described
above,
the resulting family of curves
$\tilde{\gamma}_x$ will become continuous and we will be done.

Note that some of the curves obtained in the
procedure described above could have self-intersections, however, this happens
only when they include arcs of $\alpha$ traversed twice in opposite
directions. It is easy to see that one can make all closed curves
in $\tilde G$ simple using an arbitrarily small perturbation.

\noindent {\bf Step 2.}
We will modify the homotopy $H(t,x)$ to obtain $\tilde{H}(t,x)$, a monotone homotopy with the following properties:
$\tilde{H}(*,0) = \alpha$, and $\tilde{H}(*,1)$ is contained inside the disc bounded by $H(*,1)$.
Moreover, the maximal length of curves in the new homotopy will 
increase by not more than a summand that can be made arbitrarily small. 

By analogy with Step 1, the new curves in the homotopy will be
constructed by ``pushing in'' those segments of $\beta_x$ that
lie outside the disc $D_\alpha$.  Let $D_{\beta, x}$ be the closed disc
that has $\beta_x$ as its boundary, as in the hypotheses of the lemma.  It will be a procedure that
is dual to the one in Step 1.

We will denote the curves in the new homotopy by $\tilde{\beta}_x(t)=\tilde{H}(t,x)$. In particular,
$D_\alpha \subset D_{\beta,0}$, so $\tilde{\beta}_0=\alpha$.

Now, let us describe the curve $\tilde{\beta}_1$.  If $D_{\beta_1} \subset D_\alpha$,
then we will let $\tilde{\beta}_1(t)=\beta_1(t)$.
If $D_\alpha \subset D_{\beta_1}$, then we will let $\tilde{\beta}_1=\alpha(t)$.
If $D_\alpha \cap D_\beta=\emptyset$, then we will let $\tilde{\beta}_1(t)$
be some point $\tilde{p}$, where $\tilde{p}$ is obtained as follows.
Let $\tilde{x}_0=\sup \{ x \in [0,1]$ such that 
$ D_{\beta,x} \cap D_\alpha \neq \emptyset \}$.  Let 
$\tilde{p}=D_{\beta,\tilde{x}_0} \cap D_\alpha$.

Finally, suppose that $D_\alpha \cap D_{\beta_1} \neq \emptyset$, but that one is not
a subset of the other.  In this case, $\tilde{\beta_1}$ is constructed
as follows.
Let us consider arcs of $\beta_1(t)$ that are outside
$D_\alpha$. That is, let $0=t_0<t_1<...<t_n=1$ be a subdivision of
the unit interval, such that $\beta_1(t_i)=\alpha(s_{j(i)})$ for
some $s_{j(i)} \in [0,1]$, and $\beta_1 \neq \alpha$ otherwise.  Let us consider each $\beta_1|_{[t_i, t_{i+1}]}$
that lies outside $D_\alpha$.  By the monotonicity of $H(x,t)$, it
lies in the annulus bounded by $H(0,t)$ and $G(1,t)$. The points $\alpha(s_{j(i)})$ and $\alpha(s_{j(i+1)})$
subdivide $\alpha$ into two arcs.  Exactly one of these arcs has the property that, if it is 
replaced by $\beta_1|_{[t_i,t_{i+1}]}$, then the resulting curve will be
non-contractible in the annulus.  We replace $\beta_1|_{[t_i, t_{i+1}]}$
by this arc of $\alpha$. When this is done for all arcs, we obtain $\tilde{\beta}_1$.  Clearly, 
$D_{\tilde{\beta}_1} \subset D_{\beta_1}$.

$\tilde{\beta}_x$ is constructed in a completely analogous manner for
an arbitrary $x \in [0,1]$. Note that while the length of curves in
the one-parameter family $\tilde{\beta}_x$ with $x \in [0,1]$ has not increased
compared to $\beta_x$, at this stage there can arise some discontinuities
with respect to $x$.  As in Step 1, those discontinuities can only
occur at the points where $\alpha$ and $\beta_x$ are tangent, and
only when $\beta_x$ touches $\alpha$ from the outside of $D_\alpha$, since
only arcs of $\beta_x$ outside $D_\alpha$ are to be replaced.

Note also that if $\alpha$ and $\beta_{x_0}$ intersect at only one point, 
as in Figures \ref*{monotonecase2}(a) and \ref*{monotonecase2}(b), then continuity at 
$x_0$ remains intact.

\realfig{monotonecase2}{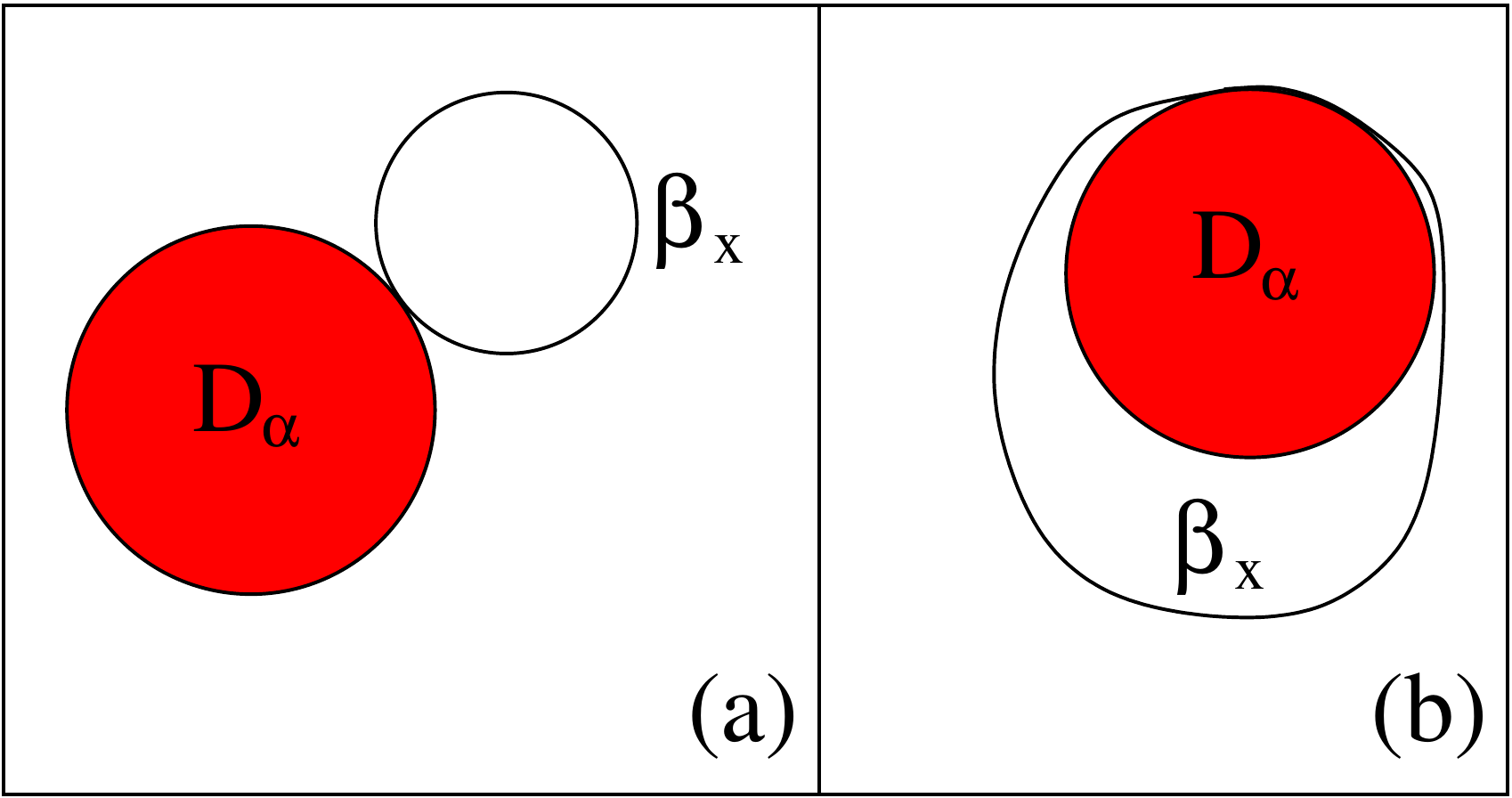}{One point of intersection}{0.7\hsize}

Thus, let us assume that there are $2k+1$ intersection points at
$0 \leq t_1 < t_2 <...<t_{2k+1} \leq 1$, $k \geq 1$.  Let $\beta(t_{j_0})=\alpha(s_{j_0})=l$ 
be the point of tangency.  Once again, let $m$ and $k$ be its neighbours with
respect to $\beta_{x_0}$.  Let us denote the arc of $\beta_{x_0}$ that connects
$m$ with $l$ by $a$ and the arc that connects $l$ with $k$ by $b$, as in 
Figures \ref*{monotoneoutside}(a) and \ref*{monotoneoutside}(b).

\realfig{monotoneoutside}{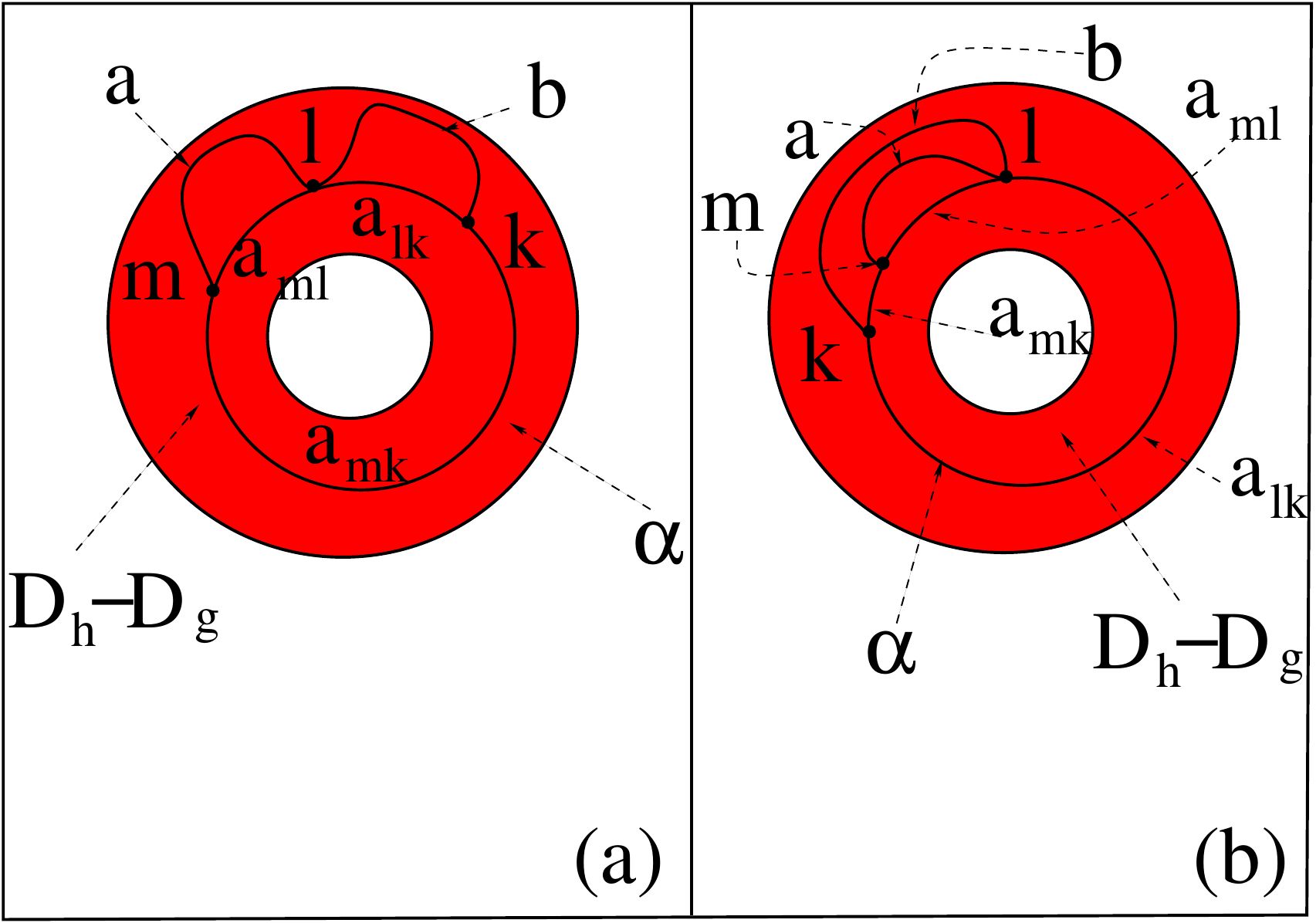}{Outside arcs}{0.7\hsize}

$m,l$, and $k$ subdivide $\alpha$ into three non-intersecting arcs
that will be denoted as $a_{ml}, a_{mk}$, and $a_{lk}$, indexed by the pair
of points that each segment connects.  Again, we are assuming that
these three arcs have disjoint interiors. Their complements in $\alpha$ will 
be denoted as $\tilde{a}_{ml}, \tilde{a}_{mk}$, and $\tilde{a}_{lk}$, respectively.

As in Case 1, a discontinuity can only arise if replacing $a$ followed by $b$ yields a different
curve than that obtained by replacing $a * b$ as a single arc.

It is easy to see that, when $\beta_{x_0}$ touches $\alpha$ from the
outside of $D_\alpha$, there are two situations to consider.

The first situation is depicted in Figure \ref*{monotoneoutside}(a).
In this case, we replace $a$ by $a_{ml}$ and $b$ by $a_{lk}$.
The other approach would be to replace the arc
$c=a*b$ by the arc $\tilde{a}_{mk}=a_{ml}*a_{lk}$.  Thus, in this
case, the two different options of replacing arcs lead to 
the same result.  In this situation, no discontinuities are formed.

The second situation is depicted in Figure \ref*{monotoneoutside}(b).
In this case, we replace $a$ by $a_{ml}$ and $b$ by
$\tilde{a}_{lk}$, while $c$ is replaced by $a_{mk}$.
While $a_{mk} \neq a_{ml}*\tilde{a}_{lk}$, we observe that  
$a_{ml}*\tilde{a}_{lk}=a_{ml}*\bar{a}_{ml}*a_{mk}$.
Thus, $a_{mk}$ is path homotopic to $a_{ml}*\tilde{a}_{lk}$
by simply contracting $a_{ml}*\bar{a}_{ml}$ along itself to 
the point $m$.  This forms a discontinuity in our homotopy $\tilde{H}$, but
if we also include this contraction (extended to the whole curve
as a part of our $1$-parameter family $\tilde{\beta}_x$), the
discontinuity will be resolved.  This is analogous to the method used to resolve discontinuities in Case 1.

\noindent {\bf Step 3.} After we have constructed
$\tilde{G}(t,x)$ and $\tilde{H}(t,x)$, we can concatenate 
$\tilde{G}$ and $\tilde{H}$ ($\tilde{G} * \tilde{H}$).  We obtain the desired homotopy $K$ by
slightly perturbing this weakly monotone homotopy to make it strictly monotone.

\end{Pf}

\section{Proof of Theorem \ref*{Theoremmain1}}

In this section, we prove Theorem \ref*{Theoremmain1}.  Since we may assume that Conjecture
\ref*{conj:hom_to_mono_gen} is true, we have a simple closed contractible curve $\gamma$ in $M$, and a strictly monotone contraction
$H(t,\tau): S^1 \times [0,1] \rightarrow M$ which covers
$\gamma$, and which consists of simple closed curves of length
no more than $L$.

Given a point $q \in \gamma$ and an $\epsilon > 0$,
we wish to construct a contraction of 
$\gamma$ through curves based at $q$ with the property that all curves are bounded in
length by
	$$ 3L + 2d + \epsilon$$
where $d$ is the diameter of the manifold.
We will also show that there is a specific point $q^\star \in \gamma$ such that there is a contraction
of $\gamma$ through curves based at $q^\star$ of length bounded by
	$$ 2L + 2d + \epsilon.$$

  Throughout this proof, we produce curves of length less than
or equal to $Q + \epsilon$ for some $Q > 0$, where $\epsilon > 0$ is chosen to be arbitrarily small.
When we combine two curves of length bounded in this way,
we simply write that the result has length bounded by $2Q + \epsilon$.  Although 
not strictly true, since we chose the original $\epsilon$ to be as small as desired, we can just
go back and choose it to be $\frac{\epsilon}{2}$, in which case the new inequality $2Q + \epsilon$
holds.  To improve readability, we do not mention this argument when it is used.

We will also be using the terms \emph{interior} and \emph{exterior} of $H_\tau$ and of $\gamma$, which we redefine here
for clarification:

\begin{Def}
\label{defn:interior_exterior}
Since $H_\tau$ is a monotone contraction, there is a disc $\mathbb{D} \subset M$ defined by the set of all points 
that are in the image of some curve in $H$.  For each point $\tau$, $H_\tau$ is simple and is contained in $\mathbb{D}$, and as such divides $M$
into two open regions.  Exactly one of these regions is entirely contained in $\mathbb{D}$.  This region is the interior
of $H_\tau$, and the other region is the exterior.  Similarly, since $H$ covers $\gamma$, $\gamma$ is contained in $\mathbb{D}$.
Since it is simple, $\gamma$ divides $M$ into $2$ regions, exactly one of which is entirely contained in $\mathbb{D}$.
This region is the interior of $\gamma$, and the other region is the exterior of $\gamma$.
\end{Def}

We will prove two lemmas which, when combined, will allow us to prove this theorem.  For each, we assume that $\epsilon > 0$ is fixed.

\begin{Lem}
\label{lem:stage_one}
There exists a point $x \in \gamma$ and a point $\tau^\star$ such that there exists a homotopy $\widetilde{H}$ from $\gamma$ to either a curve formed by slightly perturbing $H_{\tau^\star}$ or to the point $x$
through curves of length at most $2L + \epsilon$.  Additionally, $x$ lies on every curve in the homotopy $\widetilde{H}$.
\end{Lem}

Since the point $x \in \gamma$ has the aforementioned properties, $\widetilde{H}$ is a based loop homotopy.
Our second lemma takes $\widetilde{H}$ and transforms it into a contraction of $\gamma$ through curves based
at $x$ of length at most $2L + 2d + \epsilon$.

\begin{Lem}
\label{lem:stage_two}
If Lemma \ref*{lem:stage_one} does not contract $\gamma$ to $x$, then there exists a
contraction of the curve formed by slightly perturbing $H_{\tau^\star}$ through loops based at $x$ of length bounded by $2L + 2d + \epsilon$.
We denote this contraction by $K$.
\end{Lem}

We will first demonstrate how these two lemmas can be used to prove Theorem \ref*{Theoremmain1}, and will then prove each of them in turn.

\begin{Pf}{Proof of Theorem \ref*{Theoremmain1}}
Let $H$ be our original homotopy, $\widetilde{H}$ be the homotopy generated by Lemma \ref*{lem:stage_one},
and let $K$ be the homotopy generated by Lemma \ref*{lem:stage_two}.
By Lemma \ref*{lem:stage_one}, either $\widetilde{H}$ contracts $\gamma$ to the point $x$, or it homotopes $\gamma$ to a slight perturbation of $H_{\tau^\star}$.
If it contracts $\gamma$ to the point $x$, then we are done.  If it doesn't, then we have to use Lemma \ref*{lem:stage_two}.
We do this by concatenating $\widetilde{H}$ and $K$ to get a contraction of $\gamma$ through curves based at $x$ of length at most
$2L + 2d + \epsilon$, as desired.  Hence, the point $x$ is the special base point $q^\star \in \gamma$
mentioned above.  Furthermore, this will complete the proof of the theorem: if we choose any point $q \in \gamma$, then
we can build the appropriate contraction based at $q$ as follows.  Let $\alpha$ be an arc of $\gamma$ from $q$ to $x$ of length at most $\frac{L}{2}$,
and let $-\alpha$ be the same arc, but with opposite orientation.  Lastly,
let $\beta$ be the curve formed by concatenating $\alpha$ with $-\alpha$.  We can then take our contraction of $\gamma$ based at
$q^\star$, and for each curve $\gamma_\tau$ in this contraction, we replace $\gamma_\tau$ with the curve that is formed by
traversing $\alpha$, then $\gamma_\tau$, then $- \alpha$.  In this way, we produce
a homotopy from $\gamma$ to $\beta$ which is based at $q$, and which consists of curves of length at most $3L + 2d + \epsilon$.
Since $\beta$ can be contracted through loops based at $q$ of length at most $L$, this completes the proof.
\end{Pf}

We are left now with proving each of the two lemmas outlined above.

\subsection{Proof of Lemma \ref*{lem:stage_one}}

To prove Lemma \ref*{lem:stage_one}, we will adopt an approach
that will be very similar to that used by Chambers and Liokumovich in \cite{CL1}.  To begin with, we would like to
perturb the homotopy $H$ so that only finitely many non-transverse intersections between $H$
and $\gamma$ occur, and so that they do not occur concurrently.

\begin{Lem}[Perturbation Lemma]
\label{lem:perturb}
For any $\epsilon > 0$, we can perturb $H$, obtaining a new homotopy $\overline{H}$ and points
	$$ 0 = \tau_0 < \dots < \tau_n = 1 $$
such that, for all $\tau \in [\tau_i, \tau_{i+1}]$, all intersections between $H_\tau$ and $\gamma$ are transverse, except for exactly
one intersection at one point $\tau$.  The two possible interactions are shown in Figure \ref*{fig:reidemeister_move_type_B}.
$\overline{H}$ also has the following additional properties:
\begin{enumerate}
	\item	$\overline{H}$ is a contraction that covers $\gamma$.
	\item	$\overline{H}$ is monotone.
	\item	$\overline{H}$ consists of curves of length at most $L + \epsilon$.
\end{enumerate}
\end{Lem}

To prove this lemma, we use the same technique as in Proposition 2.1 from \cite{CL1}; we apply the parametric
version of Thom's Multijet Transversality Theorem to the submanifold of the $2$-fold $1$-jet bundle corresponding to curves with singularities
to show that a perturbation is possible which satisfies the above criteria.
This approach does not rule out other singular behaviour which involves self-intersections
in $\gamma$ or in $H_\tau$, however, since both of these are simple, they have no self-intersections, and so the interactions
between the two curves are limited to the isolated tangential intersections shown in Figure \ref*{fig:reidemeister_move_type_B}.
We use the term \emph{Reidemeister move} to describe this behaviour, this term
being derived from the obvious relationship between this singularity and the knot moves used in Reidemeister's Theorem.
We also note that, since $\overline{H}$ is a contraction that covers $\gamma$, $n \geq 2$.  In other
words, there must be at least $2$ Reidemeister moves, once where $\overline{H}$ transitions from a curve which lies completely in the exterior
of $\gamma$ to one which only partly lies in the exterior, and a Reidemeister move in which $H$ goes from being a curve which only partly
lies in the interior of $\gamma$, to a curve which lies either entirely in the interior of $\gamma$, or entirely in the exterior of $\gamma$.

\realfig{fig:reidemeister_move_type_B}{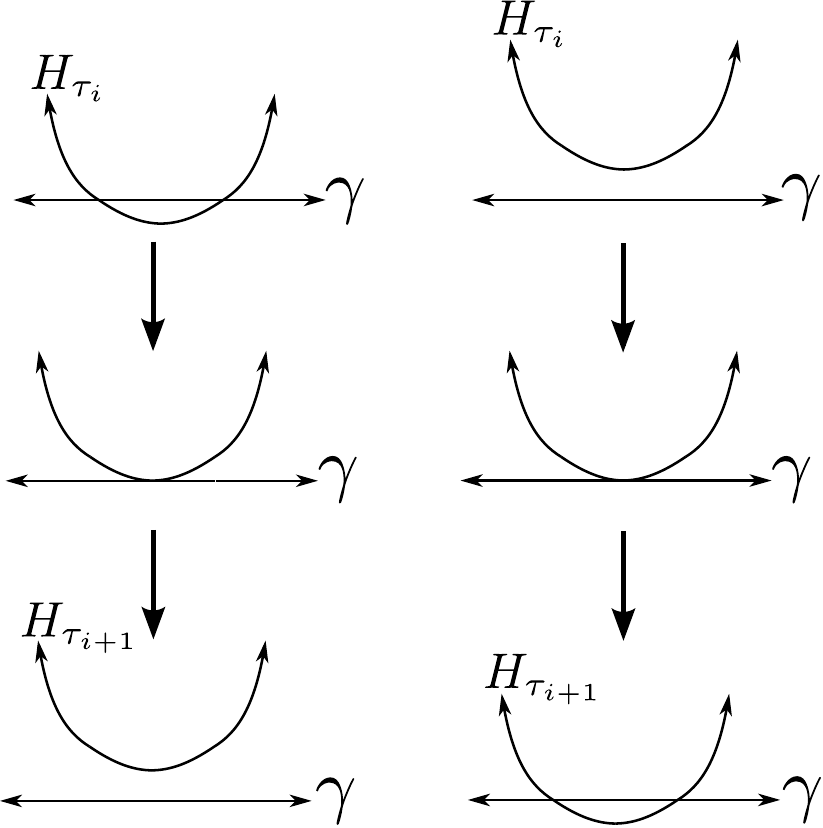}{Interactions between $\gamma$ and $H$}{1.00\textwidth}

To simplify this exposition, we will assume that $H$ has already been perturbed, and so it already has all of the properties
described in Lemma \ref*{lem:perturb}.

We now want to prove Lemma \ref*{lem:stage_one} for $H$ and $\gamma$.  We will define the point $x$ and the point $\tau^\star$, and then prove
that these points satisfy all of the required criteria.

\begin{Def}
\label{defn:x_and_t_star}
Let $x$ be the last point at which $H$ and $\gamma$ intersect, and let $\tau^\star$ be
the point at which this intersection occurs.  Note that $\tau^\star \in (\tau_{n-1}, \tau_n)$, and
$H_{\tau^\star}$ and $\gamma$ intersect tangentially at $x$.
\end{Def}

The idea to prove that Lemma \ref*{lem:stage_one} holds for these values of $x$ and $\tau^\star$ is similar to that used by
Chambers and Liokumovich in \cite{CL1}.  We construct
a certain graph $\Gamma$ where the vertices represent curves, and the edges represent homotopies between curves.  We then show that this graph
contains a certain path which represents a homotopy that easily implies the existence of the desired homotopy.

\noindent \textbf{Vertices}

We begin to construct this graph $\Gamma$ by defining its vertices.  As above, each vertex will correspond to a certain curve.
For each $i \in \{1, \dots, n - 1 \}$, consider
	$$ U_i = \gamma \cup H_{\tau_i}. $$
We will begin by identifying certain closed curves whose images lie in $U_i$.  We will then eliminate some of these 
curves based on several criteria.  For each curve that remains, we will add a vertex.  We begin by defining our large set of closed curves.
We will call these curves \emph{subcurves} at $\tau_i$.

\begin{Def}[Subcurves at $\tau_i$]
\label{defn:curves}
Choose any pair $(p_1, p_2)$ of distinct intersection points between $H_{\tau_i}$ and $\gamma$.  We can write $\gamma$
as the disjoint union of $p_1$, $p_2$, and two open segments $\rho_1$ and $\rho_2$.  Each of these segments can be used to join $p_1$ to
$p_2$.  Similarly, $H_\tau$ can be written as the disjoint union of $p_1$, $p_2$, and two open segments $\sigma_1$ and $\sigma_2$.  Each of these
segments can also be used to join $p_1$ to $p_2$.

For any piecewise smooth closed curve $\alpha$ whose image lies in $U_i$, if we can find such a pair $(p_1,p_2)$ of intersection points such that
$\alpha$ can be written as the disjoint union of $p_1$, $p_2$, $\sigma_i$ and $\rho_j$ for $i, j \in \{ 1, 2 \}$, then we say that $\alpha$
is a \emph{subcurve} at $\tau_i$.

We say that such a subcurve has endpoints $p_1$ and $p_2$, and we will denote the segment
of the curve that comes from $\gamma$ as $\rho$, and the segment that comes from $H_{\tau_i}$ as $\sigma$.  Both are open, contiguous segments
of their respective curves.
\end{Def}

Before we define which subcurves we will use to generate vertices, we will need a few definitions first.  To start, we want to define
two open, disjoint, contiguous segments of $\gamma$, which we will call $\eta_{\textrm{start}}$ and $\eta_{\textrm{end}}$.  Note that
the monotonicity of $H$ guarantees that they are disjoint.

\begin{Def}[$\eta_{\textrm{start}}$ and $\eta_{\textrm{end}}$]
\label{defn:etas}
We define the segment $\eta_{\textrm{start}}$ as the segment of $\gamma$ that is \emph{not}
contained in the closure of the interior of $H_{\tau_1}$.  Since there are exactly two intersection points between
$H_{\tau_1}$ and $\gamma$, this segment is well defined.

We define $\eta_{\textrm{end}}$ as the open segment of $\gamma$
that is contained in the interior of $H_{\tau_{n-1}}$.  Since $H_{\tau_{n-1}}$ and $\gamma$ intersect
in exactly two points, this segment is well defined.  These are shown in Figure \ref*{fig:etas}.
\end{Def}

\realfig{fig:etas}{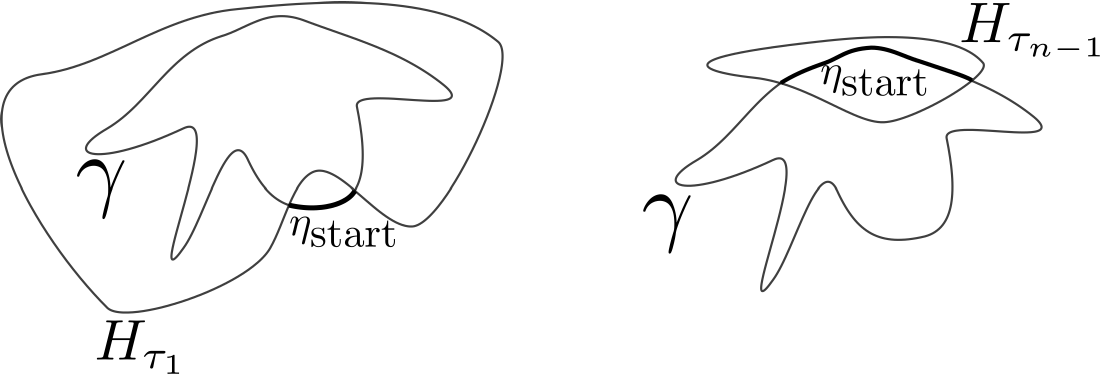}{$\eta_{\textrm{start}}$ and $\eta_{\textrm{end}}$}{1.00\textwidth}

We have a simple property of $\eta_{\textrm{start}}$ and $\eta_{\textrm{end}}$ which results from the monotonicity of $H$:
\begin{Lem}
\label{lem:etas_no_intersect}
For every point $\tau \in [\tau_1, \tau_{n-1}]$, and for any intersection point $p$ between $H_\tau$ and $\gamma$,
$p$ lies neither in $\eta_{\textrm{start}}$, nor does it lie in $\eta_{\textrm{end}}$.
\end{Lem}

We can now begin to define the set of subcurves that we will use to produce our vertices; we will define 
whether or not a subcurve \emph{respects} $\gamma$.

\begin{Def}[Respects $\gamma$]
\label{defn:respects_gamma}
We say that a subcurve $\alpha$ at $\tau_i$ \emph{respects $\gamma$} if the 
segment $\rho$ of $\alpha$ (the segment that came from $\gamma$) has the following two properties:
\begin{enumerate}
	\item	$\eta_{\textrm{start}} \cap \rho = \emptyset$
	\item	$\eta_{\textrm{end}} \subset \rho$
\end{enumerate}
\end{Def}

For every subcurve $\alpha$ at $\tau_i$ that respects $\gamma$, we give each endpoint of this curve a sign, either a $+$, or a $-$.
Let $p$ be an endpoint of $\gamma$.  Orienting $\gamma$, we can list the order in which we encounter intersection points.
Let $q$ and $r$ be the intersection points which we encounter immediately before and after $p$, which may be the same point.
Since $\gamma$ is oriented, we can also produce two contiguous segments of $\gamma$: the segment traversed from $q$ to $p$,
and the segment traversed from $p$ to $r$ (with respect to the orientation of $\gamma$).
Neither segment contains any intersection points.  Let them be $\beta_1$ and $\beta_2$.

We also see that exactly one of $\beta_1$ and $\beta_2$ must be contained in the interior of $H_{\tau_i}$ since $p$ is a transverse intersection point
of $H_{\tau_i}$ and $\gamma$.  Let this component be $\beta_j$.
Furthermore, recalling that $\rho$ is the segment of $\alpha$ that comes from $\gamma$, exactly one of $\beta_1$ and $\beta_2$ must be contained in $\rho$.
Let this component be $\beta_k$.

If $k = j$, then we assign a $+$ sign to $p$.  If not, then we assign a $-$ sign to $p$.  Note that the sign of a point does not depend on how we orient $\gamma$.

Figure \ref*{fig:intersection_sign} depicts curves $H_\tau$ and $\gamma$.  It also depicts a subcurve $\alpha$ which respects $\gamma$
and its endpoints.  Here, $\gamma$ is the same curve that appears in Figure \ref*{fig:etas}, and we assume that $\eta_{\textrm{start}}$ and
$\eta_{\textrm{end}}$ are as in this figure.  The segment $\rho$ of $\alpha$ is shown, and the signs of both endpoints are displayed as well.  Lastly,
$\sigma$ is shown with the tangent vector at each of its endpoints.
We see that the directions of these tangents
with respect to the interior of $\gamma$ do not agree with the signs of both intersection points, as per Definition \ref*{defn:valid}.  Hence, $\alpha$ is not a valid
subcurve.

\realfig{fig:intersection_sign}{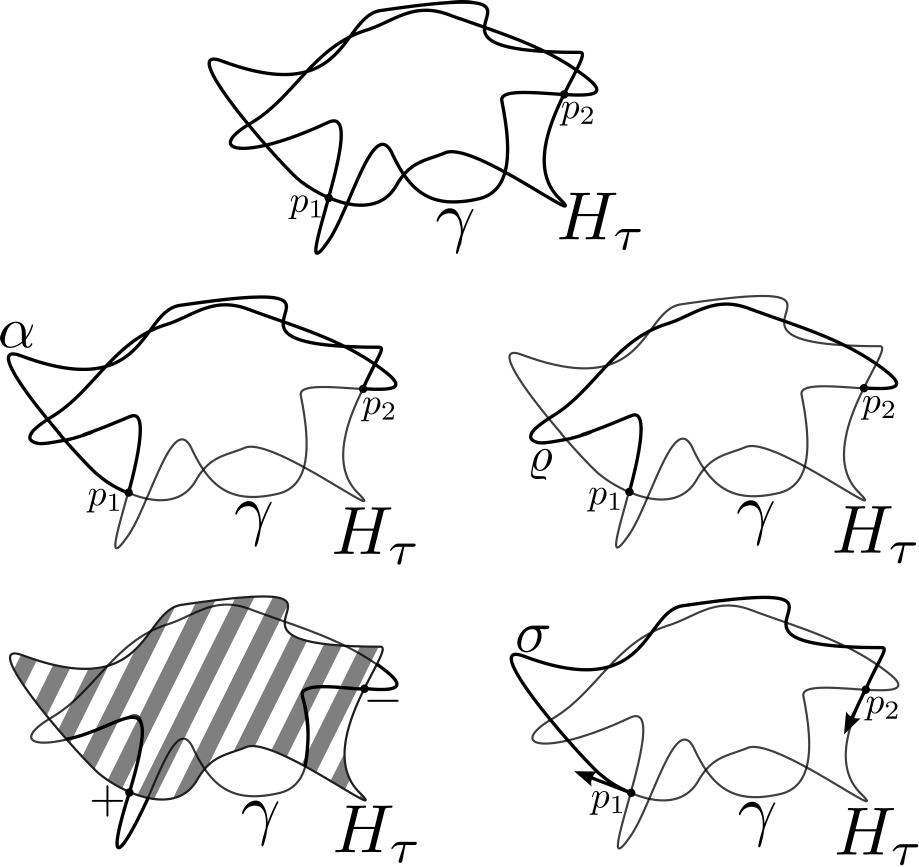}{From top to bottom, left to right: $H_\tau$ and $\gamma$, a subcurve $\alpha$ that respects $\gamma$, $\rho$,
	the signs of the endpoints of $\alpha$, and $\sigma$ with the tangent vector at each of its endpoints}{1.00\textwidth}

We can now define the set of subcurves which we want to use to produce vertices.  We call such a subcurve \emph{valid}.

\begin{Def}[Valid subcurve]
\label{defn:valid}
We say that a subcurve $\alpha$ at $\tau_i$ is \emph{valid} if it respects $\gamma$, and if the following additional properties
are true of $\sigma$.

To begin, let $p_1$ and $p_2$ be the endpoints of $\alpha$.
We can parametrize $\sigma$ so that it goes from $p_1$ to $p_2$.  Since $H_{\tau_i}$ and $\gamma$ intersect transversely at $p_1$ and at $p_2$,
we can categorize the tangent vector of $\sigma$ at $p_1$ and at $p_2$ as being into the interior of $\gamma$, or into the exterior
of $\gamma$.

We then require that, at $p_1$, the tangent of $\sigma$ points into the interior of $\gamma$ if the sign at $p_1$ is $+$,
and that it points into the exterior of $\gamma$ if it is $-$.  We also require that, at $p_2$, the tangent of $\sigma$ points into the exterior
of $\gamma$ if the sign at $p_2$ is $+$, and that the tangent points into the interior at $p_2$ if the sign is $-$.
Note that this definition is independent of the order in which we choose the endpoints of $\alpha$; a subcurve is valid with respect to one
order of endpoints if and only if it is valid with respect to the other order.
\end{Def}
For each valid subcurve at $\tau_i$ with $i \in \{ 1, \dots, n - 1 \}$, we add a vertex $v$ to the graph $\Gamma$.  We say that this vertex is generated
from $\tau_i$.  We also have a length bound
for each valid subcurve, as a result of it being composed of a segment of $\gamma$ and a segment of $H_{\tau_i}$:

\begin{Lem}[Length bound for valid subcurves]
\label{lem:length_subcurves}
For each valid subcurve $\alpha$, the length of $\alpha$ is at most
	$$ 2L + \epsilon.$$
\end{Lem}

\noindent \textbf{Edges}

We now add edges to this graph.  The idea will be that, for each $i \in \{1, \dots, n - 2 \}$, we will add a set of edges, denoted by
$E_i$.  We will specify an algorithm which takes any vertex $v$ generated from $\tau_i$ or from $\tau_{i+1}$, and produces a different vertex $w$, also
generated from $\tau_i$ or from $\tau_{i+1}$.  This algorithm is symmetric in that, if given vertex $w$, it will produce vertex $v$.
We then join each pair of vertices produced by this algorithm by an edge.  $E_i$ will be the collection of these edges.

To define this algorithm, fix a vertex $v$ in $\Gamma$ generated from $\tau_i$.
We will define the algorithm in two parts, depending on whether the resulting vertex
$w$ is generated from $\tau_{i+1}$ (a ``vertical'' edge), or if it is generated from $\tau_i$ (a ``horizontal'' edge).

Throughout the definition of this algorithm, we say that two intersections $p$ and $q$ between $H$ and $\gamma$ at $\tau_i$ are ``involved'' or ``deleted'' in the 
Reidemeister move between $\tau_i$ and $\tau_{i+1}$.  By this, we mean the following.  Let the point at which $H$ and $\gamma$ become tangent to each other be $\tau'$,
with $\tau_i < \tau' < \tau_{i+1}$.  Since all intersections between $H$ and $\gamma$ are transverse on $(\tau_i, \tau')$, we can trace the path of $p$ and
$q$ forward to $\tau'$.  When we do this, we see that $p$ gets traced to the tangential intersection at $\tau'$ (which is deleted), and
$q$ gets traced to the same intersection point.  We use the same terminology to describe intersection points between $H$ and $\gamma$ at $\tau_{i+1}$ that can be traced backwards to
the tangential intersection at $\tau'$.

\noindent \textbf{Vertical Edges}
Recall that, between $\tau_i$ and $\tau_{i+1}$, there is exactly one Reidemeister move.  This move involves two intersection points
between $H$ and $\gamma$; it either creates two intersection points, or it deletes two of them.
Let $\alpha$ be the valid subcurve at $\tau_i$ that produced the vertex $v$ and let $p_1$ and $p_2$ be the two distinct endpoints of $\alpha$.
If neither of these points is involved in the Reidemeister move, then the algorithm to find the vertex $w$
generated from $\tau_{i+1}$ is simple.  Since neither $p_1$ nor $p_2$ are deleted from $\tau_i$ to $\tau_{i+1}$, they both
follow continuous paths from $\tau_i$ to $\tau_{i+1}$.  Let $\widetilde{p_1}$ and $\widetilde{p_2}$
be the points which we reach at $\tau_{i+1}$.  We also see that we can follow $\sigma$ and $\rho$ from $\tau_i$ to $\tau_{i+1}$
in a similar fashion, arriving at $\widetilde{\sigma}$ and $\widetilde{\rho}$.  Let $\widetilde{\alpha}$ be the subcurve formed by following
$\widetilde{\sigma}$ from $\widetilde{p_1}$ to $\widetilde{p_2}$, followed by $\widetilde{\rho}$ from $\widetilde{p_2}$ back to $\widetilde{p_1}$.
If $\widetilde{\alpha}$ is a valid subcurve, then it corresponds a vertex $w$.  We will show that it is indeed valid; $w$ is the 
vertex that is produced by the algorithm.

As a result of Lemma \ref*{lem:etas_no_intersect}, $\widetilde{\rho}$ has all of the required inclusion/exclusion properties with respect to
$\eta_{\textrm{start}}$ and $\eta_{\textrm{end}}$, and so $\widetilde{\alpha}$ respects $\gamma$.  To show that it is valid, we notice 
that the sign of $p_1$ is the same as that of $\widetilde{p_1}$, and the sign of $p_2$ is the same as that of $\widetilde{p_2}$.  Let us
orient $\sigma$ from $p_1$ to $p_2$, and $\widetilde{\sigma}$ from $\widetilde{p_1}$ to $\widetilde{p_2}$.  We then have that
the direction of the tangent vector of $\sigma$ at $p_1$ with respect to the interior of $\gamma$ is the same as the direction of the tangent
vector of $\widetilde{\sigma}$ at $\widetilde{p_1}$ with respect to the interior of $\gamma$.  Similarly, the direction of the tangent vector
at $p_2$ is the same as that at $\widetilde{p_2}$.  Hence, $\widetilde{\alpha}$ is a valid subcurve, and so we are done.

If $v$ is instead generated from $\tau_{i+1}$, and neither of the endpoints of $\alpha$ are involved in the Reidemeister move
between $\tau_i$ and $\tau_{i+1}$, then we follow the exact same procedure as above, but in reverse.

\noindent \textbf{Horizontal Edges}
Again, let the vertex $v$ be generated from $\tau_i$, let $\alpha$ be the valid subcurve which corresponds to $v$, 
and let $p_1$ and $p_2$ be the endpoints of $\alpha$.  If neither $p_1$ nor $p_2$ are involved in the
Reidemeister move between $\tau_i$ and $\tau_{i+1}$, then we use the algorithm described above.
In this component of the algorithm, we determine the resulting vertex $w$ if $p_1$ or $p_2$ are involved
in the move.  Furthermore, let $\tau'$ be the point between $\tau_i$ and $\tau_{i+1}$ at which
$H_{\tau'}$ is tangent to $\gamma$.  This is the point at which the Reidemeister move ``occurs''.

We first rule out the possibility that both $p_1$ and $p_2$ are involved in the Reidemeister move:
\begin{Lem}
\label{lem:no_double_interaction}
$p_1$ and $p_2$ cannot both be deleted in the Reidemeister move between $\tau_i$ and $\tau_{i+1}$.
\end{Lem}
\begin{Pf}{Proof}
Assume that they are both involved in the 
Reidemeister move.  As in the definition of subcurves, we can break $\gamma$ into two contiguous segments, each with endpoints $p_1$ and $p_2$.
We do this by starting at $p_1$, and then by traversing $\gamma$ to $p_2$ in each of the two possible directions.  Let these two components be
$\beta_1$ and $\beta_2$.
If both $p_1$ and $p_2$ are involved in the Reidemeister move, then at least one of these two segments would have to contain no intersection points.
This is because no intersection points are deleted between $\tau_i$ and $\tau'$, and there is no way for intersection points to move through each other.
As such, until $\tau'$, the order of intersection points as we traverse $\gamma$ remains the same.  Hence, if there were intersection points
in both $\beta_1$ and $\beta_2$, then there would be no way for $p_1$ and $p_2$ to be deleted together, as there would have to be an interaction between
at least one other pair of intersection points first.  Let this intersection-free segment be $\kappa$.

We can also choose this segment $\kappa$ so that, for every $s \in \kappa$, $s$ is an intersection point between $H_{\tau_s}$ and $\gamma$
for some $\tau_s \in [\tau_i, \tau']$.

Furthermore, since $\alpha$ is a valid subcurve, it respects $\gamma$, and so we see that exactly one $\beta_j$ must 
contain $\eta_{\textrm{start}}$, and the other must contain $\eta_{\textrm{end}}$.  Hence, $\kappa$ must contain one of these curves.
By Lemma \ref*{lem:etas_no_intersect}, there are thus points in $\kappa$ that are not realized as intersection points between $\tau_i$ and $\tau_{i+1}$.
This is a contradiction, and so $p_1$ and $p_2$ cannot both be deleted in the Reidemeister move between $\tau_i$ and $\tau_{i+1}$.
\end{Pf}

Let us now move to the case where just one of $p_1$ or $p_2$ is deleted at $\tau'$.  Without loss of generality, let us assume that it is $p_1$, and let $q$ be
the other intersection point at $\tau_i$ which is deleted with $p_1$ in the Reidemeister move.  We adopt a similar approach to when we added vertical edges.
We can trace the path of $p_1$ forward until $\tau'$, and we also trace the path of $q$ forward until $\tau'$.  We notice that, since both are
deleted at $\tau'$, they merge at this point.  We can thus trace a path from $p_1$ to $q$ by first going forward to $\tau'$, tracing the path
of $p_1$ forward, and then we can go backward, tracing the path of $q$ backwards.

This path from $p_1$ to $q$ induces a homotopy from $\alpha$ to a subcurve $\widetilde{\alpha}$ at $\tau_i$ with endpoints $q$ and $p_2$.  $\widetilde{\alpha}$
is formed from the segment $\widetilde{\rho}$ of $\gamma$ and the segment $\widetilde{\sigma}$ of $H_{\tau_i}$.  The first is found by following $\rho$ forward
to $\tau'$, then by going backwards to $\tau_i$, using $q$ as an endpoint instead of $p_1$ as we go backwards; $\widetilde{\sigma}$ is found by doing the same, but with
$\sigma$.  This is shown in Figure \ref*{fig:back_and_forth}.

\realfig{fig:back_and_forth}{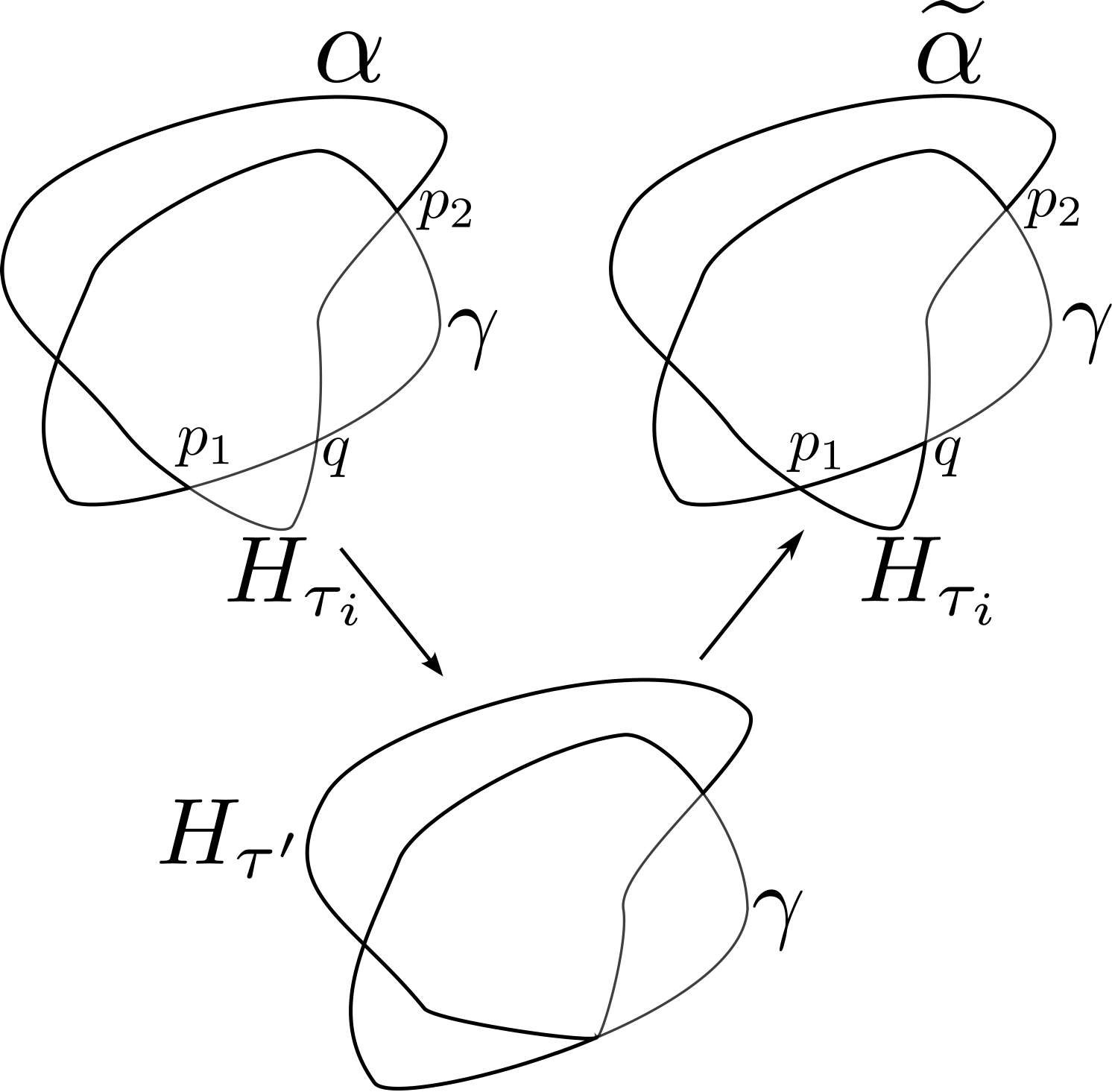}{$\alpha$ and $\widetilde{\alpha}$ if an endpoint is deleted}{1.00\textwidth}

The question, as before, is if $\widetilde{\alpha}$ is a valid subcurve.  We see that, since $\alpha$ is a valid subcurve, it respects $\gamma$, and so
$\rho$ has the proper inclusion/exclusion properties with respect to $\eta_{\textrm{start}}$ and $\eta_{\textrm{end}}$.  Lemma \ref*{lem:etas_no_intersect} then
implies that $\widetilde{\rho}$ has similar properties, and so $\widetilde{\alpha}$ respects $\gamma$.

To show that $\widetilde{\alpha}$ is valid, we must show that $\widetilde{\sigma}$ agrees with the signs of $q$ and $p_2$.  We first observe that the sign of
$q$ with respect to $\widetilde{\alpha}$ is opposite to the sign of $p_1$ with respect to $\alpha$.  On the other hand, the sign of $p_2$ remains unchanged.
If we orient $\sigma$ from $p_1$ to $p_2$, and $\widetilde{\sigma}$ from $q$ to $p_2$, then we
see that the tangent vector of $\widetilde{\alpha}$ at $p_2$ has the same direction with respect to the interior of $\gamma$ as the tangent vector of 
$\sigma$ at $p_2$, and so this endpoint meets the necessary criteria.  In terms of $q$, we see that the direction of the tangent vector of $\widetilde{\sigma}$ at $q$ with respect
to the interior of $\gamma$ is opposite to that of the tangent vector of $\sigma$ at $p_1$.  Hence, this endpoint meets the necessary criteria as well, and so $\widetilde{\alpha}$
is valid.  The rigorous proof of this is a case-by-case analysis on the segment of $H_{\tau_i} \cup \gamma$ around $q$ and $p_1$.  The cases are formed by
considering all possible interiors of $H_{\tau_i}$, and all possible arcs $\rho$.  This is shown in Figure \ref*{fig:intersection_fix}.
Since $\widetilde{\alpha}$ is valid, it corresponds to a vertex $w$, which is the desired vertex.

\realfig{fig:intersection_fix}{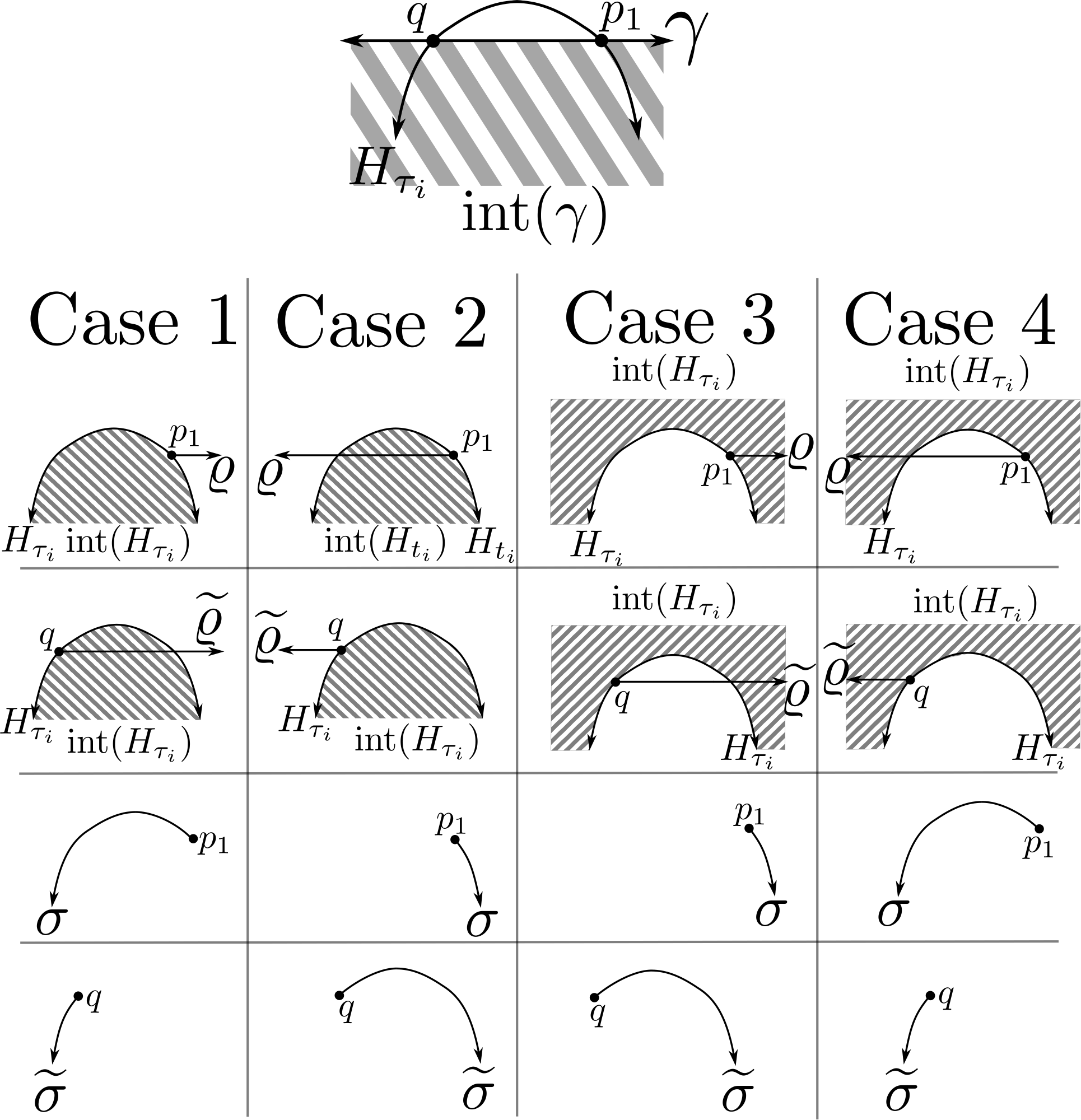}{In order of rows, top to bottom: intersections $q$ and $p_1$ with the interior of $\gamma$ shaded, $\rho$ with the interior of $H_{\tau_i}$ shaded,
	   $\widetilde{\rho}$ with the interior of $H_{\tau_i}$ shaded, $\sigma$, and $\widetilde{\sigma}$}{1.00\textwidth}

If $v$ is generated from $\tau_{i+1}$, then we follow the above procedure, but in reverse.  That is, if the Reidemeister move
between $\tau_i$ and $\tau_{i+1}$ creates two intersection points of which one is an endpoint of $\alpha$, then we follow the above steps
to produce a vertex $w$.  Note that for reasons analogous to those presented in the proof of Lemma \ref*{lem:no_double_interaction}, both endpoints of $\alpha$ cannot be created
by the Reidemeister move between $\tau_i$ and $\tau_{i+1}$.

Before we complete the proof of Lemma \ref*{lem:stage_one}, we prove some important properties of $\Gamma$:

\begin{Lem}[Properties of $\Gamma$]
\label{lem:gamma_properties}
The graph $\Gamma$ has the following properties:
\begin{enumerate}
	\item	For each set of edges $E_i$, $i \in \{1, \dots, n - 2 \}$, and for each vertex $v$ generated from $\tau_i$ or $\tau_{i+1}$,
		$v$ is the endpoint of exactly one edge in $E_i$.
	\item	All vertices generated from $\tau_1$ and all vertices generated from $\tau_{n-1}$ have degree 1;
		all other vertices have degree $2$.
	\item	There is exactly one vertex generated from $\tau_1$, and one vertex generated from $\tau_{n-1}$,
		and they correspond to the curves shown in Figure \ref*{fig:endpoint_curves}.
	\item	If two vertices are joined by an edge, then there is a homotopy of closed curves between
		the curves corresponding to the vertices through closed curves of length at most $2L + \epsilon$.
		Furthermore, all of these curves contain $\eta_{\textrm{end}}$.
\end{enumerate}
\end{Lem}
\begin{Pf}{Proof}
The first statement is a result of the fact that the algorithm used to add edges takes any vertex $v$ generated from $\tau_i$ or from 
$\tau_{i+1}$ and produces a vertex $w$, $v \neq w$.  Since we use this algorithm to add edges, there is an edge between $v$ and $w$.  Additionally, it is easy to check
that this algorithm is symmetric in that the vertex $w$ will produce the vertex $v$.  Hence, each vertex is the endpoint of exactly one
edge in $E_i$.

The second statement results from the fact that, for each vertex $v$ that is generated from $\tau_1$ or from $\tau_{n-1}$, $v$ is the endpoint
of exactly one edge from $E_1$ or $E_{n-1}$, respectively, and there is no other set $E_j$ which contains an edge that has $v$ as an endpoint.
Thus, the degree is $1$.  For any vertex $v$ generated from $\tau_i$ with $i \in \{ 1, \dots, n-1 \}$, $v$ is the endpoint of an edge from
$E_i$, and is also the endpoint of an edge from $E_{i+1}$.  Hence, it has degree $2$.

The third statement follows from looking at the set of all valid subcurves at $\tau_1$ and $\tau_{n-1}$.  At each of these points, there are exactly two
intersections between $H$ and $\gamma$, and so it is a simple exercise to look at each of the four subcurves and to show that the only ones
that are valid are the ones depicted in Figure \ref*{fig:endpoint_curves}.

The last statement comes from examining how we add edges.  In all of the cases, we are tracing two intersection points back or forth, and
keeping track of one segment of $H_\tau$ that connects these $2$ points and one segment of $\gamma$ that also connects these two points, which generates a continuous homotopy.  Since
both $\gamma$ and $H_\tau$ are bounded in length by $L + \epsilon$, taking a segment of one and joining it with a segment from the other has length
at most
	$$ 2L + \epsilon.$$
The fact that they all contain $\eta_{\textrm{end}}$ is a result of two observations.  First, all subcurves at any $\tau_i$ respect $\gamma$, and so contain $\eta_{\textrm{end}}$.
Second, as a result of Lemma \ref*{lem:etas_no_intersect}, no intersections between $H$ and $\gamma$ lie in $\eta_{\textrm{end}}$ for any $\tau \in [\tau_1, \tau_{n-1}]$.
\end{Pf}

\realfig{fig:endpoint_curves}{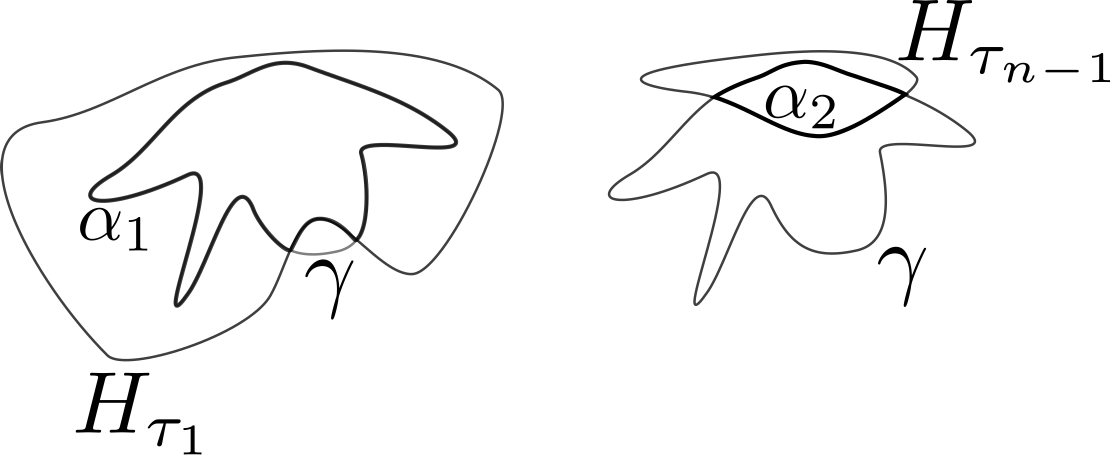}{The curve that corresponds to the only vertex generated from $\tau_1$, and the curve that corresponds to the only vertex generated from $\tau_{n-1}$}{1.00\textwidth}

We can now prove Lemma \ref*{lem:stage_one}.

\begin{Pf}{Proof of Lemma \ref*{lem:stage_one}}
From Lemma \ref*{lem:gamma_properties}, we have that there is only $1$ vertex $v$ generated from $\tau_1$, and one vertex $w$ generated from $\tau_{n-1}$.  Additionally, they have degree $1$, and all other
vertices in $\Gamma$ have degree $2$.  As a result, we have that there is a path in $\Gamma$ from $v$ to $w$.
Let $\alpha_1$ and $\alpha_2$ be the subcurves that correspond to $v$ and $w$, respectively.
Due to the property of $\Gamma$ that edges represent homotopies over closed curves of length at most $2L + \epsilon$,
there is thus a homotopy from $\alpha_1$ to $\alpha_2$ over such curves.  Furthermore, every curve in this homotopy contains
$\eta_{\textrm{end}}$.

We now observe that $\gamma$ is homotopic to $\alpha_1$ over curves of length at most $2L + \epsilon$, and so we can homotope
$\gamma$ to $\alpha_2$ over closed curves with the same length bound.  All of these curves also contain $\eta_{\textrm{end}}$.

The rest of the proof depends on whether $H$ contracts $\gamma$ to a point inside $\gamma$ or outside $\gamma$.  If it contracts
$\gamma$ to a point outside $\gamma$, then we see that $\alpha_2$ can be contracted to the point $x \in \eta_{\textrm{end}}$ on $\gamma$
through curves that contain $x$.
Since $\eta_{\textrm{end}}$ is contained in all curves in this homotopy up to $\alpha_2$, we conclude that $x$
is contained in every curve in this entire contraction.

If $H$ contracts $\gamma$ to a point inside $\gamma$, then recalling that $\tau^\star$ is the last point at which $H$ intersects $\gamma$, and 
$x$ is the point of tangential intersection at $\tau^\star$, we can homotope $\alpha_2$ to $H_{\tau^\star}$ through curves containing $x$ and 
which are bounded in length by $2L + \epsilon$.  Since $\gamma$ can also be homotoped to $\alpha_2$ through such curves,
this gives us a desirable homotopy from $\gamma$ to $H_{\tau^\star}$.  This completes the proof.
\end{Pf}

Finally, we illustrate this process using an explicit homotopy.  This is shown in Figure \ref*{fig:example}.

\realfig{fig:example}{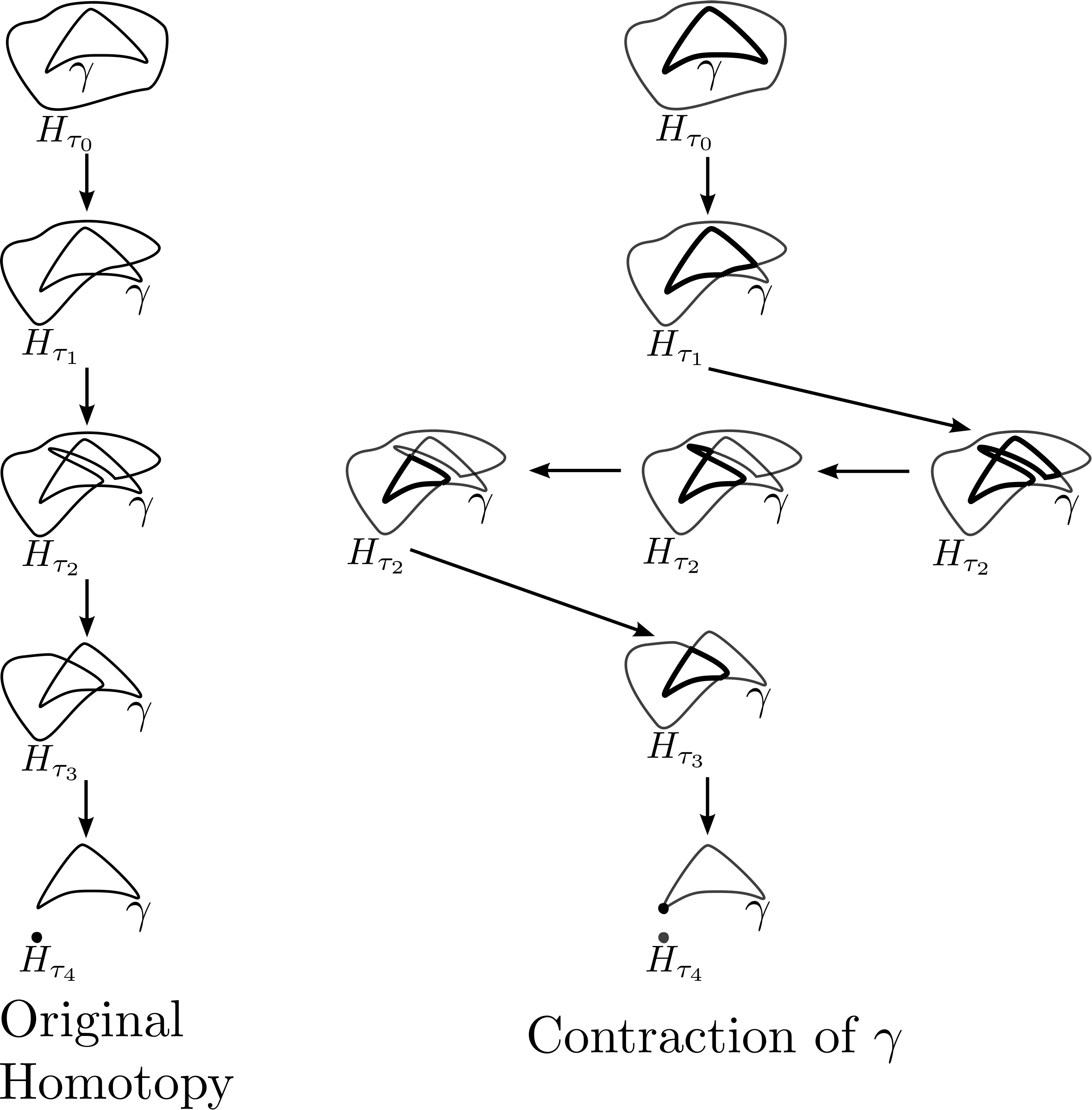}{A homotopy that covers $\gamma$ and the resulting contraction of $\gamma$}{1.00\textwidth}


\subsection{Proof of Lemma \ref*{lem:stage_two}}

We now prove Lemma \ref*{lem:stage_two}.  Given a curve $H_{\tau^\star}$ and a point $x \in \gamma \cap H_{\tau^\star}$, we want to show that we can contract $H_{\tau^\star}$ through curves based at $x$,
and of length at most $2L + 2d + \epsilon$.  The idea here is to employ
a method similar to that used in this article to produce a contraction of the boundary of a Riemannian disc from a monotone contraction of that boundary.
To do this, let $c$ be the point that $H$ contracts $H_{\tau^\star}$ to.  Join $x$ to $c$ via a minimal geodesic, and let $y$ be the last point of intersection between this geodesic
and $H_{\tau^\star}$.  This is depicted in Figure \ref*{fig:curve_pieces}.

\realfig{fig:curve_pieces}{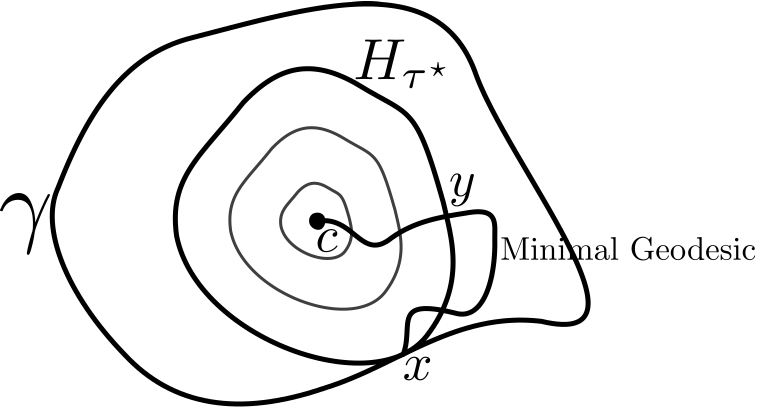}{$H$ and $\gamma$ as per the hypotheses of Lemma \ref*{lem:stage_two}}{1.00\textwidth}

Our homotopy now works as follows.  One should refer to Figure \ref*{fig:example_homotopy} for a visual reference.
Let $\beta$ be the segment of the minimal geodesic that connects $y$ to $c$ entirely in the interior of $H_{\tau^\star}$.  Let the length of $\beta$ be $B$; we of course
have that $B \leq d$, where $d$ is the diameter of the manifold.  Let $\alpha$ be a segment of $H_{\tau^\star}$ that connects $x$ to $y$ which is of length at most $\frac{L + \epsilon}{2}$.
We now produce our contraction of $H_{\tau^\star}$ in three parts.
The first part is a homotopy from $H_{\tau^\star}$ to the curve formed by traversing $\alpha$ from $x$ to $y$, following by traversing the entirety of $H_{\tau^\star}$ from $y$ to $y$, and then by traversing $- \alpha$
from $y$ back to $x$.  This homotopy consists of curves bounded in length by $2L + \epsilon$.  Let us call this curve $\eta$.

\realfig{fig:example_homotopy}{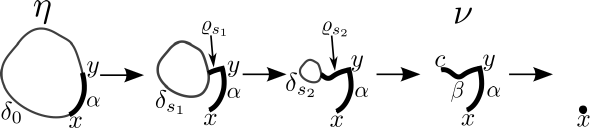}{Building the contraction}{1.00\textwidth}

The second step is a homotopy from $\eta$ to the curve formed by traversing $\alpha$ from $x$ to $y$, then $\beta$ from $y$ to $c$, then $- \beta$ from $c$ back to $y$, then $- \alpha$ from $y$ back to $x$.
Let us call this curve $\nu$.  This homotopy, $P$, is defined on the interval $[0, B]$, where (as above) $B$ is the length of $\beta$.  For each
$s \in [0,B]$, let $\rho_s$ be the segment of $\beta$ from $y$ which has a length of $s$.  Since $H$ is monotone, there is exactly one curve $\delta$ corresponding to a curve in the homotopy
$H$ which has the property that it goes through $y$ if $s = 0$, and that it goes through the endpoint of $\rho_s$ which is not $y$ if $s > 0$.  Now, we define $P(s)$ to be the curve formed by
traversing $\alpha$, then $\rho_s$, then $\delta$, then $- \rho_s$, then $- \alpha$.  Since $H$ is monotone, this produces a continuous homotopy of piecewise smooth simple curves
of length at most $2L + 2d + \epsilon$.

The third step is that we homotope $\nu$ to the point $x$ by contracting it in the obvious way; since it is a curve traversed forward from $x$ to $y$ to $c$ and then backward from $c$ to $y$ to $x$,
it is obvious how to do this without exceeding a length bound of $L + 2d + \epsilon$.

By concatenating the above homotopies, we get a homotopy of closed curves with the desired properties, completing the proof.

\bigskip

\noindent {\bf Acknowledgements:} This work was supported in part by an NSERC Discovery Grant (Rotman), by NSERC Postgraduate and Postdoctoral Scholarships (Chambers), and
by an Ontario Graduate Scholarship (Chambers). This
paper was partially written during the authors' visit of the 
Max-Planck Institute for Mathematics in Bonn.  The authors would like to
thank the Max-Planck Institute for its kind hospitality.

This work is based on Chambers' doctoral work.  The authors would like to thank
Arnaud de Mesmay, Erin W. Chambers, and Tim Ophelders for finding an error in an earlier version of this article.



 


\begin{tabbing}
\hspace*{7.5cm}\=\kill
Gregory R. Chambers                 \> Regina Rotman\\
Department of Mathematics           \> Department of Mathematics\\
University of Chicago               \> University of Toronto\\
Chicago, Illinois, 60637            \> Toronto, Ontario M5S 2E4\\
USA                                 \> Canada\\
e-mail: chambers@math.uchicago.edu  \> rina@math.utoronto.ca\\
\end{tabbing}


\end{document}